\documentclass[10pt]{article}
\usepackage{a4wide}
\usepackage{theorem}
\usepackage{amssymb}
\newcommand\al{\alpha}
\newcommand\be{\beta}
\newcommand\ga{\gamma}
\newcommand\de{\delta}

\newcommand\ro{\varrho}
\newcommand\si{\sigma}
\newcommand\ta{\tau}

\newcommand\om{\omega}

\newcommand\Ga{\Gamma}

\newcommand\Si{\Sigma}

\newcommand\Om{\Omega}

\newcommand\Hom{\mbox{\rm Hom}}

\newcommand\id{{\rm id}}

\newcommand\lbra{[\hspace{-2pt}[}
\newcommand\rbra{]\hspace{-2pt}]}

\def\NN{\mathbb N}

\def\QQ{\mathbb Q}

\def\ZZ{\mathbb Z}

  \newcounter{numero}[section]
  \newcounter{aux1}
  \newcounter{aux2}
  \renewcommand{\thenumero}{\thesection.\arabic{numero}}
  \newenvironment{corollary}{
  \dimen255=\parindent \parindent=0in
  \refstepcounter{numero}{\vspace{.3cm}\bf Corollary \thenumero .}\
  \it}{\rm \parindent=\dimen255\par}

  \newenvironment{theorem}{
  \dimen255=\parindent \parindent=0in
  \refstepcounter{numero}{\vspace{.3cm}\bf Theorem \thenumero .}\
  \it}{\rm \parindent=\dimen255\par}

  \newenvironment{lemma}{
  \dimen255=\parindent \parindent=0in
  \refstepcounter{numero} {\vspace{.3cm}\bf Lemma \thenumero .}\
  \it}{\rm \parindent=\dimen255\par}

  \newenvironment{proposition}{
  \dimen255=\parindent \parindent=0in
  \refstepcounter{numero} {\vspace{.3cm}\bf Proposition \thenumero .}\ \it
  }{\rm \parindent=\dimen255\par}

  \newenvironment{proof}{
  {\vspace{.3cm}\noindent\sc{Proof:}}\
  }{$\fbox{}$ \par}

\title{Vassiliev invariants for braids on surfaces}
\author{ Juan Gonz\'alez-Meneses \quad  Luis Paris}
\date{May, 2000}

\begin{document}
\maketitle
\sloppy

\begin{abstract}
 We show that Vassiliev invariants separate braids on a closed 
oriented surface, and we exhibit an universal Vassiliev invariant for 
these braids in terms of chord diagrams labeled by elements of the 
fundamental group of the considered surface.
\footnote{\hspace{-.65cm} Keywords: Braid - Surface -  
Vassiliev Invariant - Finite Type Invariant.}
\footnote{\hspace{-.65cm} {\em Mathematics Subject Classification:} 
 Primary: 20F36. Secondary: 57M27, 57N05.}
\footnote{\hspace{-.65cm} First author partially supported by 
DGESIC-PB97-0723 and 
by the european network TMR Sing. Eq. Diff. et Feuill.}
\end{abstract}

\section{Definitions and statements}\label{secdef}

\subsection{Introduction}
  
  Vassiliev knot invariants were introduced by V. A. Vassiliev 
(\cite{vassiliev1}, \cite{vassiliev2}; see also \cite{birman2},
\cite{bar-natan}), and they have been generalized 
to several other knot-like objects, such as links, braids, tangles, 
string links, knotted graphs, etc. The purpose of this paper is to 
consider Vassiliev invariants of braids on surfaces, and to extend
some well known results of Vassiliev invariants of Artin braids
to the case of braids on surfaces.

 Our study of Vassiliev invariants is inspired by Papadima's work 
\cite{papadima} on Vassiliev invariants for Artin braids with values in 
$\ZZ$. However, the presence of the fundamental group of the surface
varies substantially the analysis. Anyway, the Vassiliev theory for 
braids on surfaces, exposed in this paper, appears to be a natural 
generalization of the corresponding theory for Artin braids.

\noindent {\bf Aknowledgement.} We are grateful to \c Stefan Papadima 
for stimulating conversations and suggestions which were the starting 
point of this work.

\subsection{Braids and singular braids on surfaces}

  Throughout this paper $M$ will denote a closed, orientable surface of genus 
$g\geq 1$, and ${\cal P}=\{P_1,\ldots,P_n \}$  a set of $n$ distinct 
points in $M$. Define a {\em $n$-braid based at ${\cal P}$} to be a 
collection $b = (b_1,\ldots,b_n)$ of disjoint smooth paths in 
$M\times [0,1]$, called {\em strings} of $b$, such that the $i$-th string 
$b_i$ runs monotonically in $t\in [0,1]$ from the point $(P_i,0)$ to some 
point $(P_j,1)$, $\; P_j\in {\cal P}$.

  An {\em isotopy} in this context is a deformation through braids
(which fixes the ends).  Multiplication of braids is defined by 
concatenation, generalizing the construction of the fundamental group.
The isotopy classes of braids with this multiplication form the group 
$B_n(M,{\cal P})$, called {\em braid group with $n$ strings on $M$ based
at ${\cal P}$}. Note that the group $B_n(M,{\cal P})$ does not depend,
up to isomorphism, on the set ${\cal P}$ of points but only on the
cardinality $n=|{\cal P}|$. So we may write $B_n(M)$ in place of 
$B_n(M,{\cal P})$.

 In the same way as Artin braid groups have been extended to singular 
braid monoids (\cite{birman2},\cite{baez}), one can extend the braid group 
$B_n(M)$ to $SB_n(M)$, the {\em monoid of singular braids with 
$n$ strings on $M$}. The strings of a singular braid 
are now allowed to intersect transversely, but only in finitely many double 
points, called {\em singular points}.

  As with braids, isotopy is a deformation through singular braids 
(which fixes the ends), and multiplication is by concatenation. Note
that the isotopy classes of singular braids form a monoid and not a group:
the singular braids with one or more singular points being non-invertible.

\subsection{Vassiliev invariants and Vassiliev filtration}
 
  An {\em invariant} of braids on $M$ with values in an abelian group $A$,
is a set-mapping $v:\: B_n(M)\rightarrow A$. 
Like for knots and Artin braids, one can extend $v$ to singular braids by 
using the recursive rule

\vspace{.3cm}
\begin{figure}[ht]
\centerline{\begin{picture}(0,0)%
\includegraphics{vresol75.pstex}%
\end{picture}%
\setlength{\unitlength}{2960sp}%
\begingroup\makeatletter\ifx\SetFigFont\undefined
\def\x#1#2#3#4#5#6#7\relax{\def\x{#1#2#3#4#5#6}}%
\expandafter\x\fmtname xxxxxx\relax \def\y{splain}%
\ifx\x\y   
\gdef\SetFigFont#1#2#3{%
  \ifnum #1<17\tiny\else \ifnum #1<20\small\else
  \ifnum #1<24\normalsize\else \ifnum #1<29\large\else
  \ifnum #1<34\Large\else \ifnum #1<41\LARGE\else
     \huge\fi\fi\fi\fi\fi\fi
  \csname #3\endcsname}%
\else
\gdef\SetFigFont#1#2#3{\begingroup
  \count@#1\relax \ifnum 25<\count@\count@25\fi
  \def\x{\endgroup\@setsize\SetFigFont{#2pt}}%
  \expandafter\x
    \csname \romannumeral\the\count@ pt\expandafter\endcsname
    \csname @\romannumeral\the\count@ pt\endcsname
  \csname #3\endcsname}%
\fi
\fi\endgroup
\begin{picture}(4554,736)(766,-744)
\put(766,-376){\makebox(0,0)[lb]{\smash{\SetFigFont{14}{16.8}{rm}$v$}}}
\put(4006,-376){\makebox(0,0)[lb]{\smash{\SetFigFont{14}{16.8}{rm}$- \; v$}}}
\put(2161,-376){\makebox(0,0)[lb]{\smash{\SetFigFont{14}{16.8}{rm}$= \; v$}}}
\end{picture}
}
\label{vresol}
\end {figure}

\vspace{.3cm}
  The picture on the left hand side represents a small neighborhood
of a singular point in a singular braid. Those on the right hand side
represent the braids which are obtained from the previous one by 
{\em resolution} of that singular point. That is, we modify the first
braid inside the neighborhood of the singular point, in a {\em positive} 
and a {\em negative} way, to obtain two singular braids having one less
singular point.

  Remark that this is well defined since $M$ is orientable. In a 
non-orientable case, one also has two different modifications, but it 
would not be possible to differentiate the positive one from the
negative one.

 Let $d$ be an integer. A {\em Vassiliev invariant of type $d$} is
an invariant $v$ such that $v(b)=0$ for every singular braid $b$ with
more than $d$ singular points.

  There is an equivalent definition of a Vassiliev invariant, in terms
of the so-called {\em Vassiliev filtration}. First, consider the group ring
$\ZZ[B_n(M)]$. We can define a map 
$$
  \eta:\: SB_n(M)\longrightarrow  \ZZ[B_n(M)]
$$
which ``resolves'' all the singular points of a given braid, with the
corresponding signs. That is,

\vspace{.3cm}
\begin{figure}[ht]
\centerline{\begin{picture}(0,0)%
\includegraphics{eta75.pstex}%
\end{picture}%
\setlength{\unitlength}{2960sp}%
\begingroup\makeatletter\ifx\SetFigFont\undefined
\def\x#1#2#3#4#5#6#7\relax{\def\x{#1#2#3#4#5#6}}%
\expandafter\x\fmtname xxxxxx\relax \def\y{splain}%
\ifx\x\y   
\gdef\SetFigFont#1#2#3{%
  \ifnum #1<17\tiny\else \ifnum #1<20\small\else
  \ifnum #1<24\normalsize\else \ifnum #1<29\large\else
  \ifnum #1<34\Large\else \ifnum #1<41\LARGE\else
     \huge\fi\fi\fi\fi\fi\fi
  \csname #3\endcsname}%
\else
\gdef\SetFigFont#1#2#3{\begingroup
  \count@#1\relax \ifnum 25<\count@\count@25\fi
  \def\x{\endgroup\@setsize\SetFigFont{#2pt}}%
  \expandafter\x
    \csname \romannumeral\the\count@ pt\expandafter\endcsname
    \csname @\romannumeral\the\count@ pt\endcsname
  \csname #3\endcsname}%
\fi
\fi\endgroup
\begin{picture}(4679,742)(744,-713)
\put(4411,-376){\makebox(0,0)[lb]{\smash{\SetFigFont{14}{16.8}{rm}$-$}}}
\put(2341,-151){\makebox(0,0)[lb]{\smash{\SetFigFont{11}{13.2}{rm}$\eta$}}}
\end{picture}
}
\label{eta}
\end {figure}

\vspace{.3cm}
  This map is a well defined multiplicative morphism. Remark that a singular 
braid with $d$ singular points is mapped to an alternate sum of $2^d$ 
non-singular braids, each one having coefficient $+1$ or $-1$ depending on 
the sign of its corresponding resolutions.

  Let $S_dB_n(M)$ denote the set of isotopy classes of 
singular braids with $d$ singular points. We denote by $V_d$ the 
$\ZZ$-module generated by $\eta(S_dB_n(M))$.
One can easily verify that $V_d$ is a (two-sided) ideal of $\ZZ[B_n(M)]$
and that we have the inclusions $V_{d+1}\subset V_d$ and 
$V_{d_1}V_{d_2}\subset V_{d_1+d_2}$, for all $d_1,d_2,d\in \NN$. 
We have then obtained a filtration
$$
  \ZZ[B_n(M)]= V_0 \supset V_1 \supset V_2 \supset \cdots,
$$
which is called the {\em Vassiliev filtration} of $\ZZ[B_n(M)]$.

  The definition of a Vassiliev invariant in terms of the Vassiliev 
filtration is as follows. One can extend any invariant 
$v: B_n(M) \rightarrow A$ by linearity to
a morphism of $\ZZ$-modules $v:\: \ZZ[B_n(M)]\rightarrow A$.
Note that the previous extension of $v$ to singular braids can also 
be expressed by $v(b)= v(\eta(b))$,
for $b\in SB_n(M)$.   Then, $v$ is a Vassiliev invariant of type 
$d$ if and only if it vanishes on $V_{d+1}$. Therefore, the set of 
Vassiliev invariants of type $d$ with values in $A$ is equal to 
$$
\Hom_{\ZZ}(\ZZ[B_n(M)]/V_{d+1}\:,\:A).
$$

\subsection{Statements}

  We have two goals in this paper. The first one is to show that Vassiliev
invariants separate braids on surfaces, that is, to prove the following.

\vspace{.3cm}
\begin{theorem}\label{separa}
Given two non-equivalent braids $b$ and $c$ on $M$, there exists an integer 
$N\geq 1$ and a Vassiliev invariant $v_N$ of type $N$ such that 
$v_N(b)\neq v_N(c)$. Moreover, $v_N$ can be chosen to take values in 
$\ZZ$.
\end{theorem}

\vspace{.3cm} 
  This result is known to hold for Artin braids (\cite{bar-natan2},
\cite{kohno}, \cite{papadima}), but it is still a conjecture for knots.
Actually, Theorem~\ref{separa} is a corollary of the following theorem.

\vspace{.3cm}
\begin{theorem}\label{prosepara}
 Let $\{V_d\}_{d=1}^{\infty}$ be the Vassiliev filtration of 
$\ZZ[B_n(M)]$. One has:
\begin{enumerate}
 \item $\bigcap_{d=0}^{\infty} V_d = \{ 0 \}$,

 \item $V_d/V_{d+1} \;$ is a free $\ZZ$-module for all $d\geq 0$.
\end{enumerate}
\end{theorem}

\vspace{.3cm} 
  Indeed, if Theorem~\ref{prosepara} holds, then given two 
non-equivalent braids $b,c\in B_n(M)$, there exists an integer $N$
such that $b-c\notin V_{N+1}$. Then we can take $v_N$ to be the canonical 
projection from $\ZZ[B_n(M)]$ to $\ZZ[B_n(M)]/V_{N+1}$.
In addition, if $V_d/V_{d+1}$ is a free $\ZZ$-module for all $d$, then 
$$
  \ZZ[B_n(M)]/V_{N+1}\simeq (\ZZ[B_n(M)]/V_1) \oplus (V_1/V_2) \oplus \cdots
   \oplus (V_N/V_{N+1})
$$
is also a free $\ZZ$-module, so
we can obviously compose the above projection with a map from 
$\ZZ[B_n(M)]/V_{N+1}$ to $\ZZ$, in such a way that the image of 
$b-c$ is non-zero. Therefore, our first goal will be achieved 
by proving Theorem~\ref{prosepara}.

\begin{figure}[ht]
\centerline{\begin{picture}(0,0)%
\includegraphics{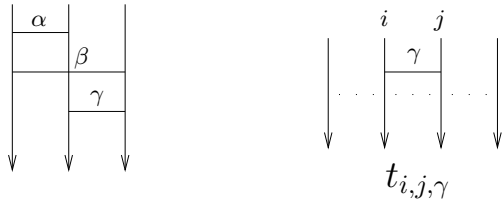}%
\end{picture}%
\setlength{\unitlength}{2960sp}%
\begingroup\makeatletter\ifx\SetFigFont\undefined
\def\x#1#2#3#4#5#6#7\relax{\def\x{#1#2#3#4#5#6}}%
\expandafter\x\fmtname xxxxxx\relax \def\y{splain}%
\ifx\x\y   
\gdef\SetFigFont#1#2#3{%
  \ifnum #1<17\tiny\else \ifnum #1<20\small\else
  \ifnum #1<24\normalsize\else \ifnum #1<29\large\else
  \ifnum #1<34\Large\else \ifnum #1<41\LARGE\else
     \huge\fi\fi\fi\fi\fi\fi
  \csname #3\endcsname}%
\else
\gdef\SetFigFont#1#2#3{\begingroup
  \count@#1\relax \ifnum 25<\count@\count@25\fi
  \def\x{\endgroup\@setsize\SetFigFont{#2pt}}%
  \expandafter\x
    \csname \romannumeral\the\count@ pt\expandafter\endcsname
    \csname @\romannumeral\the\count@ pt\endcsname
  \csname #3\endcsname}%
\fi
\fi\endgroup
\begin{picture}(3894,1596)(889,-1132)
\put(1576,-331){\makebox(0,0)[lb]{\smash{\SetFigFont{9}{10.8}{rm}$\ga$}}}
\put(4006,299){\makebox(0,0)[lb]{\smash{\SetFigFont{9}{10.8}{rm}$i$}}}
\put(4456,299){\makebox(0,0)[lb]{\smash{\SetFigFont{9}{10.8}{rm}$j$}}}
\put(4231, 29){\makebox(0,0)[lb]{\smash{\SetFigFont{9}{10.8}{rm}$\ga$}}}
\put(4051,-1051){\makebox(0,0)[lb]{\smash{\SetFigFont{14}{16.8}{rm}$t_{i,j,\ga}$}}}
\put(1081,299){\makebox(0,0)[lb]{\smash{\SetFigFont{9}{10.8}{rm}$\al$}}}
\put(1441,-16){\makebox(0,0)[lb]{\smash{\SetFigFont{9}{10.8}{rm}$\be$}}}
\end{picture}
}
\caption{A $M$-labeled chord diagram and the generator $t_{i,j,\ga}$.}
\label{chordM}
\end {figure}

  Our second goal is to define a universal Vassiliev invariant for $B_n(M)$
which generalizes the notion of chord diagrams for Artin braids. 
Recall that a {\em chord diagram} is a diagram made of $n$ vertical lines 
and of a finite number of horizontal segments, called chords, connecting
the lines. A {\em $M$-labeled chord diagram} is a chord diagram such that
each chord is labeled by an element of $\pi_1(M)$ (see Figure~\ref{chordM}).
Note that the set of $M$-labeled chord diagrams is equipped with a
multiplication defined by concatenation. The free $\ZZ$-module
generated by the chord diagrams is a $\ZZ$-algebra which can be
identified with $\ZZ[t_{i,j,\ga}]$, the free non-commutative $\ZZ$-algebra
freely generated by the $t_{i,j,\ga}$, where $i,j\in\{1,\ldots,n \}$, 
$i\neq j$, $\:\ga\in \pi_1(M)$, and where $t_{i,j,\ga}=t_{j,i,\ga^{-1}}$ 
(see Figure~\ref{chordM}).

\vspace{.3cm}

 We denote by ${\cal A}_n$ the quotient $\ZZ$-algebra obtained from
$\ZZ[t_{i,j,\ga}]$ by imposing the relations

\begin{itemize}
 \item  $[t_{i,j,\ga}\; ,\; t_{k,l,\de}]=0,$ \quad for all
distinct $ i,j,k,l\in \{1,\ldots, n\}$ and all $\ga, \de \in \pi_1(M)$,

 \item $[t_{i,j,\ga}\; ,\; t_{j,k,\de}+t_{i,k,(\ga\de)}]=0$, \quad 
for all distinct $ i,j,k \in \{1,\ldots, n\}$ and all $\ga, \de \in \pi_1(M)$,
\end{itemize}
and we denote by $\widehat{\cal A}_n$ its natural completion.

 Note that the symmetric group $\Si_n$ acts on $\pi_1(M)^n$ by 
permuting coordinates, so we can consider the induced semi-direct product
$H_n=\pi_1(M)^n\rtimes\Si_n$. In addition, it is straightforward to show 
that $H_n$ acts on $\widehat{\cal A}_n$, defining the semi-direct 
product  $\widehat{\cal A}_n\rtimes \ZZ[H_n]$. The action is defined 
by the following relations.

\begin{itemize}

 \item $\si \; t_{i,j,\ga}\; \si^{-1} = t_{\si(i),\si(j),\ga}$, \quad
   for all $\si\in \Si_n$,

 \item $\mu(k)\; t_{i,j,\ga}\; \mu(k)^{-1} = t_{i,j,\ga}$, \quad
   for all $\mu\in \pi_1(M)$ and all $k\neq i,j$,
 
 \item $\mu(i)\; t_{i,j,\ga} \;\mu(i)^{-1} = t_{i,j,(\mu\ga)}$, \quad
   for all $\mu \in \pi_1(M)$,

\end{itemize}
where, $\mu(i)= (1,\ldots,1,\stackrel{(i)}{\mu},1,\ldots,1)\in \pi_1(M)^n$.
Note that one also has the following relation:
$$
 \mu(j)\; t_{i,j,\ga}\; \mu(j)^{-1} = \mu(j) \; t_{j,i,\ga^{-1}} \;
\mu(j)^{-1} =  t_{j,i,(\mu\ga^{-1})}= t_{i,j,(\ga\mu^{-1})}.
$$

  The $\ZZ$-algebra $\widehat{\cal A}_n\rtimes \ZZ[H_n]$ carries the 
filtration induced by that of $\widehat{\cal A}_n$, so its associated graded 
algebra is ${\cal A}_n\rtimes \ZZ[H_n]$. We also have 
$\;\mbox{gr}_V\ZZ[B_n(M)]=\bigoplus_{d=0}^{\infty}(V_d/V_{d+1})$.
Our second main result will be:

\vspace{.3cm}
\begin{theorem}\label{teorisom}
There exists a homomorphism of $\ZZ$-modules 
$\;u: \ZZ[B_n(M)] \rightarrow 
\widehat{\cal A}_n \rtimes \ZZ[H_n]$ such that the corresponding graded
map 
$$
   \mbox{\em gr}u: \: \mbox{\em gr}_V\ZZ[B_n(M)] \longrightarrow 
  {\cal A}_n \rtimes \ZZ[H_n]
$$
is an isomorphism of graded $\ZZ$-algebras. 
\end{theorem}

\vspace{.6cm} 
 We end this section by showing why $u$ is called a {\em universal Vassiliev
invariant} for $B_n(M)$.

\vspace{.3cm}
\begin{corollary}
  Every Vassiliev invariant of $B_n(M)$ factors through $u$ in a unique way.
\end{corollary}

\vspace{.3cm}
\begin{proof}
  By Theorem~\ref{prosepara}, we know that $\ZZ[B_n(M)]/V_{N+1}$ is 
a free $\ZZ$-module for all $N\geq 0$, hence, $\ZZ[B_n(M)]\simeq 
(\ZZ[B_n(M)]/V_{N+1}) \oplus V_{N+1}$. Recall that 
$$
 \mbox{gr}_V\ZZ[B_n(M)]= \bigoplus_{d=0}^{\infty}{(V_d/V_{d+1})}
\simeq \left(\ZZ[B_n(M)]/V_{N+1}\right)\oplus 
\left(\bigoplus_{d>N}{(V_d/V_{d+1})}\right).
$$
Now, since gr$u$ is an isomorphism, we conclude that, for all $N\geq 0$, 
$\;{\cal A}_n^{(\leq N)} \rtimes \ZZ[H_n]$ is also a free $\ZZ$-module.
Therefore, one has:
$$
  \widehat{\cal A}_n\rtimes \ZZ[H_n]\simeq 
  ({\cal A}_n^{(\leq N)} \rtimes \ZZ[H_n])
 \oplus (\widehat{\cal A}_n^{(> N)} \rtimes \ZZ[H_n]),
$$
and
$$
  {\cal A}_n\rtimes \ZZ[H_n]\simeq ({\cal A}_n^{(\leq N)} \rtimes \ZZ[H_n])
 \oplus ({\cal A}_n^{(> N)} \rtimes \ZZ[H_n]).
$$

  Every Vassiliev invariant $v\in \Hom_{\ZZ}(\ZZ[B_n(M)]/V_{N+1}\:,\:A)$
can then be seen as a linear map from $\mbox{gr}_V\ZZ[B_n(M)]$ to $A$ which 
vanishes on $\bigoplus_{d>N}{(V_d/V_{d+1})}$. Via gr$u$, this means
that $v$ is a linear map from ${\cal A}_n \rtimes \ZZ[H_n]$ to $A$, which 
vanishes on ${\cal A}_n^{(> N)} \rtimes \ZZ[H_n]$. Therefore, if
$v$ is a Vassiliev invariant of type $N$, it can be lifted in a unique way 
to a linear map $\widehat{v}:\widehat{\cal A}_n\rtimes \ZZ[H_n] 
\rightarrow A$, which verifies $v= \widehat{v}\circ u$.
\end{proof}

\section{Vassiliev invariants separate braids}

Our strategy for proving Theorem~\ref{prosepara} is the following.
In a first subsection, we introduce some ideal $J$ of $\ZZ[B_n(M)]$ 
given by its generators and we prove that $V_d=J^d$ for all $d\geq 0$.
In a second subsection, we consider an exact sequence
$1\rightarrow K_n \rightarrow B_n(M) \rightarrow H_n \rightarrow 1$,
and we prove that $J^d$ is equal in some sense  to 
$I(K_n)^d\otimes \ZZ[H_n]$, where $I(K_n)$ denotes the augmentation ideal 
of $K_n$. 
 In a third subsection, we prove that $K_n$ can be expressed as an iterated
semi-direct product of free groups (of infinite rank). Finally,
in the fourth subsection, we use the results of the previous ones to
prove Theorem~\ref{prosepara}.

\subsection{The Vassiliev filtration coincides with the $J$-adic filtration}
\label{coincides}

 The aim of this subsection is to introduce an ideal $J$ of $\ZZ[B_n(M)]$
defined by its generators, and to show that the Vassiliev filtration
coincides with the $J$-adic filtration (i.e. $V_d=J^d$ for all $d\in \NN$).

  We begin by explaining our ``visualization'' of (singular) braids,
and by exhibiting generators for $SB_n(M)$. 

\begin{figure}[ht]
\centerline{\begin{picture}(0,0)%
\includegraphics{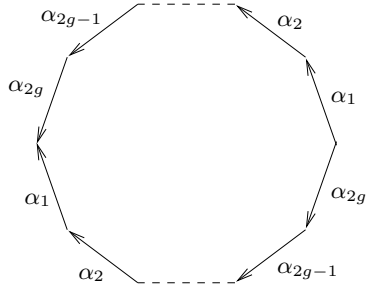}%
\end{picture}%
\setlength{\unitlength}{2960sp}%
\begingroup\makeatletter\ifx\SetFigFont\undefined
\def\x#1#2#3#4#5#6#7\relax{\def\x{#1#2#3#4#5#6}}%
\expandafter\x\fmtname xxxxxx\relax \def\y{splain}%
\ifx\x\y   
\gdef\SetFigFont#1#2#3{%
  \ifnum #1<17\tiny\else \ifnum #1<20\small\else
  \ifnum #1<24\normalsize\else \ifnum #1<29\large\else
  \ifnum #1<34\Large\else \ifnum #1<41\LARGE\else
     \huge\fi\fi\fi\fi\fi\fi
  \csname #3\endcsname}%
\else
\gdef\SetFigFont#1#2#3{\begingroup
  \count@#1\relax \ifnum 25<\count@\count@25\fi
  \def\x{\endgroup\@setsize\SetFigFont{#2pt}}%
  \expandafter\x
    \csname \romannumeral\the\count@ pt\expandafter\endcsname
    \csname @\romannumeral\the\count@ pt\endcsname
  \csname #3\endcsname}%
\fi
\fi\endgroup
\begin{picture}(2700,2383)(496,-2063)
\put(3196,-511){\makebox(0,0)[lb]{\smash{\SetFigFont{9}{10.8}{rm}$\al_1$}}}
\put(496,-421){\makebox(0,0)[lb]{\smash{\SetFigFont{9}{10.8}{rm}$\al_{2g}$}}}
\put(2746,-1951){\makebox(0,0)[lb]{\smash{\SetFigFont{9}{10.8}{rm}$\al_{2g-1}$}}}
\put(631,-1366){\makebox(0,0)[lb]{\smash{\SetFigFont{9}{10.8}{rm}$\al_1$}}}
\put(1081,-1996){\makebox(0,0)[lb]{\smash{\SetFigFont{9}{10.8}{rm}$\al_2$}}}
\put(766,164){\makebox(0,0)[lb]{\smash{\SetFigFont{9}{10.8}{rm}$\al_{2g-1}$}}}
\put(3196,-1321){\makebox(0,0)[lb]{\smash{\SetFigFont{9}{10.8}{rm}$\al_{2g}$}}}
\put(2746,119){\makebox(0,0)[lb]{\smash{\SetFigFont{9}{10.8}{rm}$\al_2$}}}
\end{picture}
}
\caption{A representation of the surface $M$.}
\label{polygon}
\end {figure}

  We represent the surface $M$ as a polygon of $4g$ sides which are 
identified in the way of Figure~\ref{polygon}.

\begin{figure}[ht]
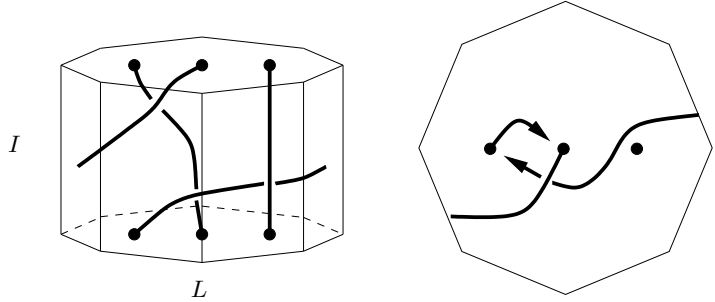

\centerline{\begin{picture}(0,0)%
\includegraphics{viewpoint175.pstex}%
\end{picture}%
\setlength{\unitlength}{2960sp}%
\begingroup\makeatletter\ifx\SetFigFont\undefined
\def\x#1#2#3#4#5#6#7\relax{\def\x{#1#2#3#4#5#6}}%
\expandafter\x\fmtname xxxxxx\relax \def\y{splain}%
\ifx\x\y   
\gdef\SetFigFont#1#2#3{%
  \ifnum #1<17\tiny\else \ifnum #1<20\small\else
  \ifnum #1<24\normalsize\else \ifnum #1<29\large\else
  \ifnum #1<34\Large\else \ifnum #1<41\LARGE\else
     \huge\fi\fi\fi\fi\fi\fi
  \csname #3\endcsname}%
\else
\gdef\SetFigFont#1#2#3{\begingroup
  \count@#1\relax \ifnum 25<\count@\count@25\fi
  \def\x{\endgroup\@setsize\SetFigFont{#2pt}}%
  \expandafter\x
    \csname \romannumeral\the\count@ pt\expandafter\endcsname
    \csname @\romannumeral\the\count@ pt\endcsname
  \csname #3\endcsname}%
\fi
\fi\endgroup
\begin{picture}(2667,2082)(271,-1906)
\put(1801,-1906){\makebox(0,0)[lb]{\smash{\SetFigFont{9}{10.8}{rm}$L$}}}
\put(271,-691){\makebox(0,0)[lb]{\smash{\SetFigFont{9}{10.8}{rm}$I$}}}
\end{picture}
\quad \quad \quad 
  \begin{picture}(0,0)%
\includegraphics{viewpoint275.pstex}%
\end{picture}%
\setlength{\unitlength}{2960sp}%
\begingroup\makeatletter\ifx\SetFigFont\undefined
\def\x#1#2#3#4#5#6#7\relax{\def\x{#1#2#3#4#5#6}}%
\expandafter\x\fmtname xxxxxx\relax \def\y{splain}%
\ifx\x\y   
\gdef\SetFigFont#1#2#3{%
  \ifnum #1<17\tiny\else \ifnum #1<20\small\else
  \ifnum #1<24\normalsize\else \ifnum #1<29\large\else
  \ifnum #1<34\Large\else \ifnum #1<41\LARGE\else
     \huge\fi\fi\fi\fi\fi\fi
  \csname #3\endcsname}%
\else
\gdef\SetFigFont#1#2#3{\begingroup
  \count@#1\relax \ifnum 25<\count@\count@25\fi
  \def\x{\endgroup\@setsize\SetFigFont{#2pt}}%
  \expandafter\x
    \csname \romannumeral\the\count@ pt\expandafter\endcsname
    \csname @\romannumeral\the\count@ pt\endcsname
  \csname #3\endcsname}%
\fi
\fi\endgroup
\begin{picture}(2364,2364)(635,-2136)
\end{picture}
}
\caption{A braid with $3$ strings on a surface of genus $2$: two 
different viewpoints.}
\label{viewpoint}
\end {figure}

  We draw braids over $M$ in this polygon, as if we looked at the
cylinder $M\times [0,1]$ from above, that is, we project the braid over
$M\times\{0\}$. Like for the planar representations of knots, we 
see over and under-crossings, and we can always move our braid 
via a convenient isotopy to avoid triple crossing points in the 
projection. See Figure~\ref{viewpoint} for an example.

\begin{figure}[ht]
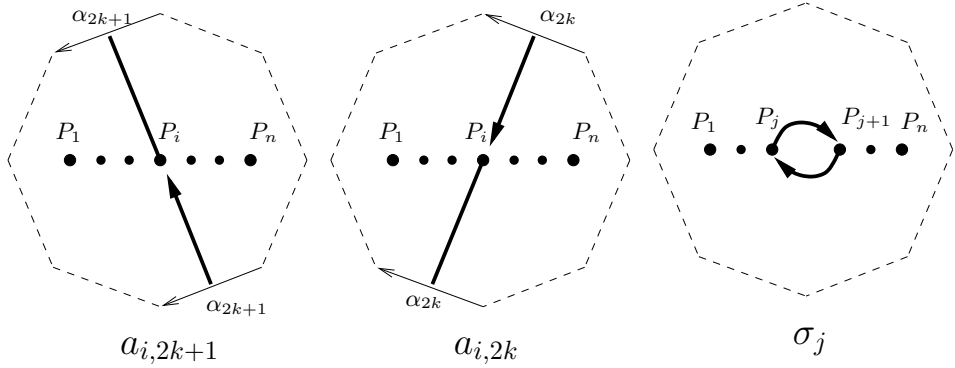

\centerline{\begin{picture}(0,0)%
\includegraphics{air75.pstex}%
\end{picture}%
\setlength{\unitlength}{2960sp}%
\begingroup\makeatletter\ifx\SetFigFont\undefined
\def\x#1#2#3#4#5#6#7\relax{\def\x{#1#2#3#4#5#6}}%
\expandafter\x\fmtname xxxxxx\relax \def\y{splain}%
\ifx\x\y   
\gdef\SetFigFont#1#2#3{%
  \ifnum #1<17\tiny\else \ifnum #1<20\small\else
  \ifnum #1<24\normalsize\else \ifnum #1<29\large\else
  \ifnum #1<34\Large\else \ifnum #1<41\LARGE\else
     \huge\fi\fi\fi\fi\fi\fi
  \csname #3\endcsname}%
\else
\gdef\SetFigFont#1#2#3{\begingroup
  \count@#1\relax \ifnum 25<\count@\count@25\fi
  \def\x{\endgroup\@setsize\SetFigFont{#2pt}}%
  \expandafter\x
    \csname \romannumeral\the\count@ pt\expandafter\endcsname
    \csname @\romannumeral\the\count@ pt\endcsname
  \csname #3\endcsname}%
\fi
\fi\endgroup
\begin{picture}(2409,2985)(619,-2491)
\put(1171,299){\makebox(0,0)[lb]{\smash{\SetFigFont{9}{10.8}{rm}$\al_{2k+1}$}}}
\put(1891,-691){\makebox(0,0)[lb]{\smash{\SetFigFont{9}{10.8}{rm}$P_{i}$}}}
\put(1036,-691){\makebox(0,0)[lb]{\smash{\SetFigFont{9}{10.8}{rm}$P_{1}$}}}
\put(2656,-691){\makebox(0,0)[lb]{\smash{\SetFigFont{9}{10.8}{rm}$P_{n}$}}}
\put(2296,-2131){\makebox(0,0)[lb]{\smash{\SetFigFont{9}{10.8}{rm}$\al_{2k+1}$}}}
\put(1576,-2491){\makebox(0,0)[lb]{\smash{\SetFigFont{14}{16.8}{rm}$a_{i,2k+1}$}}}
\end{picture}
\quad \begin{picture}(0,0)%
\includegraphics{ais75.pstex}%
\end{picture}%
\setlength{\unitlength}{2960sp}%
\begingroup\makeatletter\ifx\SetFigFont\undefined
\def\x#1#2#3#4#5#6#7\relax{\def\x{#1#2#3#4#5#6}}%
\expandafter\x\fmtname xxxxxx\relax \def\y{splain}%
\ifx\x\y   
\gdef\SetFigFont#1#2#3{%
  \ifnum #1<17\tiny\else \ifnum #1<20\small\else
  \ifnum #1<24\normalsize\else \ifnum #1<29\large\else
  \ifnum #1<34\Large\else \ifnum #1<41\LARGE\else
     \huge\fi\fi\fi\fi\fi\fi
  \csname #3\endcsname}%
\else
\gdef\SetFigFont#1#2#3{\begingroup
  \count@#1\relax \ifnum 25<\count@\count@25\fi
  \def\x{\endgroup\@setsize\SetFigFont{#2pt}}%
  \expandafter\x
    \csname \romannumeral\the\count@ pt\expandafter\endcsname
    \csname @\romannumeral\the\count@ pt\endcsname
  \csname #3\endcsname}%
\fi
\fi\endgroup
\begin{picture}(2409,2985)(619,-2491)
\put(1666,-2491){\makebox(0,0)[lb]{\smash{\SetFigFont{14}{16.8}{rm}$a_{i,2k}$}}}
\put(1036,-691){\makebox(0,0)[lb]{\smash{\SetFigFont{9}{10.8}{rm}$P_{1}$}}}
\put(2656,-691){\makebox(0,0)[lb]{\smash{\SetFigFont{9}{10.8}{rm}$P_{n}$}}}
\put(1756,-691){\makebox(0,0)[lb]{\smash{\SetFigFont{9}{10.8}{rm}$P_{i}$}}}
\put(2386,299){\makebox(0,0)[lb]{\smash{\SetFigFont{9}{10.8}{rm}$\al_{2k}$}}}
\put(1261,-2086){\makebox(0,0)[lb]{\smash{\SetFigFont{9}{10.8}{rm}$\al_{2k}$}}}
\end{picture}
\quad 
  \begin{picture}(0,0)%
\includegraphics{sij75.pstex}%
\end{picture}%
\setlength{\unitlength}{2960sp}%
\begingroup\makeatletter\ifx\SetFigFont\undefined
\def\x#1#2#3#4#5#6#7\relax{\def\x{#1#2#3#4#5#6}}%
\expandafter\x\fmtname xxxxxx\relax \def\y{splain}%
\ifx\x\y   
\gdef\SetFigFont#1#2#3{%
  \ifnum #1<17\tiny\else \ifnum #1<20\small\else
  \ifnum #1<24\normalsize\else \ifnum #1<29\large\else
  \ifnum #1<34\Large\else \ifnum #1<41\LARGE\else
     \huge\fi\fi\fi\fi\fi\fi
  \csname #3\endcsname}%
\else
\gdef\SetFigFont#1#2#3{\begingroup
  \count@#1\relax \ifnum 25<\count@\count@25\fi
  \def\x{\endgroup\@setsize\SetFigFont{#2pt}}%
  \expandafter\x
    \csname \romannumeral\the\count@ pt\expandafter\endcsname
    \csname @\romannumeral\the\count@ pt\endcsname
  \csname #3\endcsname}%
\fi
\fi\endgroup
\begin{picture}(2409,2793)(619,-2572)
\put(2206,-646){\makebox(0,0)[lb]{\smash{\SetFigFont{9}{10.8}{rm}$P_{j+1}$}}}
\put(2701,-691){\makebox(0,0)[lb]{\smash{\SetFigFont{9}{10.8}{rm}$P_{n}$}}}
\put(946,-691){\makebox(0,0)[lb]{\smash{\SetFigFont{9}{10.8}{rm}$P_{1}$}}}
\put(1801,-2491){\makebox(0,0)[lb]{\smash{\SetFigFont{14}{16.8}{rm}$\si_{j}$}}}
\put(1486,-646){\makebox(0,0)[lb]{\smash{\SetFigFont{9}{10.8}{rm}$P_{j}$}}}
\end{picture}
}
\caption{Generators for $B_n{M}$.}
\label{generators1}
\end {figure}

 Now, for every $i\in\{1,\ldots,n\}$ and every $r\in \{1,\ldots,2g\}$,  
we define the braid  $a_{i,r}$ as follows. All the strings of $a_{i,r}$
are trivial except the $i$-th one which goes through the $r$-th
wall in the way of Figure~\ref{generators1}. It goes upwards if
$r$ is odd and it goes downwards if $r$ is even.

 We also define, for all $j=1,\ldots, n-1$, the braid $\si_j$ as follows.
All the strings of $\si_j$ are trivial except the $j$-th one and the
$(j+1)$-th one. The $j$-th string goes from $(P_j,0)$ to $(P_{j+1},1)$
and the $(j+1)$-th string goes from $(P_{j+1},0)$ to $(P_j,1)$ according
to Figure~\ref{generators1}. Note that $\si_1,\ldots,\si_{n-1}$ are the 
classical generators of the braid group $B_n$ of the disc.

  It is easy to show that $\{a_{i,r};\;\; i=1,\ldots,n,\: r=1,\ldots, 2g\}
\cup \{\si_1,\ldots, \si_{n-1} \}$ is a generator set for $B_n(M)$.
Actually, there is no need to include $a_{i,r}$ if $i\geq 2$, but it is 
better for our purposes. One can find in~\cite{juan} a presentation for
$B_n(M)$ which involves these generators.

\begin{figure}[ht]
\centerline{\begin{picture}(0,0)%
\includegraphics{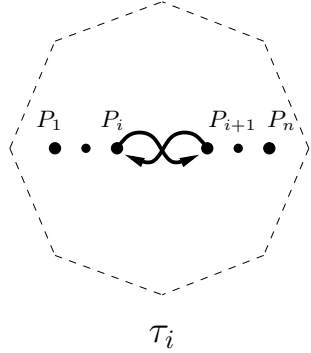}%
\end{picture}%
\setlength{\unitlength}{2960sp}%
\begingroup\makeatletter\ifx\SetFigFont\undefined
\def\x#1#2#3#4#5#6#7\relax{\def\x{#1#2#3#4#5#6}}%
\expandafter\x\fmtname xxxxxx\relax \def\y{splain}%
\ifx\x\y   
\gdef\SetFigFont#1#2#3{%
  \ifnum #1<17\tiny\else \ifnum #1<20\small\else
  \ifnum #1<24\normalsize\else \ifnum #1<29\large\else
  \ifnum #1<34\Large\else \ifnum #1<41\LARGE\else
     \huge\fi\fi\fi\fi\fi\fi
  \csname #3\endcsname}%
\else
\gdef\SetFigFont#1#2#3{\begingroup
  \count@#1\relax \ifnum 25<\count@\count@25\fi
  \def\x{\endgroup\@setsize\SetFigFont{#2pt}}%
  \expandafter\x
    \csname \romannumeral\the\count@ pt\expandafter\endcsname
    \csname @\romannumeral\the\count@ pt\endcsname
  \csname #3\endcsname}%
\fi
\fi\endgroup
\begin{picture}(2409,2712)(619,-2491)
\put(2791,-691){\makebox(0,0)[lb]{\smash{\SetFigFont{9}{10.8}{rm}$P_{n}$}}}
\put(1801,-2491){\makebox(0,0)[lb]{\smash{\SetFigFont{14}{16.8}{rm}$\tau_{i}$}}}
\put(856,-691){\makebox(0,0)[lb]{\smash{\SetFigFont{9}{10.8}{rm}$P_{1}$}}}
\put(1396,-691){\makebox(0,0)[lb]{\smash{\SetFigFont{9}{10.8}{rm}$P_{i}$}}}
\put(2296,-691){\makebox(0,0)[lb]{\smash{\SetFigFont{9}{10.8}{rm}$P_{i+1}$}}}
\end{picture}
}
\caption{The singular braid $\tau_i$.}
\label{taui}
\end {figure}

  For every $i=1,\ldots,n-1$, we define the singular braid 
$\tau_i\in S_1B_n(M)$ as in Figure~\ref{taui}. This singular braid has
a unique singular point on which intersect the $i$-th string and the
$(i+1)$-th string.  The $i$-th string goes from $(P_i,0)$ to 
$(P_{i+1},1)$ and the $(i+1)$-th string goes from $(P_{i+1},0)$
to $(P_i,1)$. The other strings are trivial.

 By a suitable isotopy, any singular braid $b\in S_kB_n(M)$ can be
written in the form
$$
    b= c_1 \ta_{j_1} c_2 \ta_{j_2} \cdots c_k \ta_{j_k} c_{k+1},
$$
where $c_i\in B_n(M)$. So the following set generates $SB_n(M)$ (as a 
monoid). 
$$
   \{a_{i,r}^{\pm 1};\;\; i=1,\ldots,n,\: r=1,\ldots, 2g  \}
   \cup \{\si_1^{\pm 1},\ldots,\si_{n-1}^{\pm 1} \}
   \cup \{\tau_1,\ldots, \tau_{n-1} \}.
$$

 Now, the morphism $\eta:\: SB_n(M) \rightarrow \ZZ[B_n(M)]$ sends
$\si_i^{\pm 1}$ to $\si_i^{\pm 1}$, $a_{i,r}^{\pm 1}$ to $a_{i,r}^{\pm 1}$ 
and $\tau_i$ to $\si_i - \si_i^{-1}$. 
Recall that $V_d$ denotes the $\ZZ$-module of
$\ZZ[B_n(M)]$ generated by $\eta(S_dB_n(M))$. From the above 
considerations, it immediately follows:

\vspace{.3cm}
\begin{proposition}
 Let $J$ be the two-sided ideal of $\ZZ[B_n(M)]$ generated by
$\{ \si_i-\si_i^{-1}; \; \; i=1,\ldots, n-1 \}$.
Then $V_d=J^d$ for all $d\in \NN$.
\end{proposition}

\subsection{From $J^d$ to $I(K_n)^d$}\label{secIJ}
  
 Recall that $H_n$ denotes the semi-direct product 
$\pi_1(M)^n\rtimes \Si_n$. We define a homomorphism
$\varphi:\:  B_n(M)  \rightarrow  H_n$ as follows. We fix a disc
$D$ embedded in $M$ which contains ${\cal P}$ and, for all 
$i,j\in\{1,\ldots,n\}$, a path $\al_{i,j}$ in $D$ going from $P_i$ 
to $P_j$. Pick a braid $b=(b_1,\ldots,b_n)$, 
$\:b_i:\:[0,1]\rightarrow M\times[0,1]$, based at ${\cal P}$.
Let $s\in \Si_n$ be the permutation induced by $b$.  
Let $\overline{b}_i:\: [0,1]\rightarrow M$ be the projection of $b_i$ 
on the first coordinate, and let $\mu_i$ be the loop based at $P_i$
defined by $\mu_i=\overline{b}_i\:\al_{s(i),i}$. Then we set
$$
    \varphi(b) = (\mu_1,\ldots,\mu_n) s \in \pi_1(M)^n\rtimes \Si_n
			= H_n.
$$

One can easily verify that $\varphi:\:  B_n(M)  \rightarrow  H_n$ is a  
well defined homomorphism, and that its definition depends
on the choice of $D$ but not on the choice of the paths $\al_{i,j}$.

  Let $K_n$ denote the kernel of $\varphi$.
It is a classical matter that a set-section  $\si:\: H_n\rightarrow B_n(M)$
of $\varphi$ determines a $\ZZ$-isomorphism
$\Phi:\: \ZZ[B_n(M)]  \rightarrow \ZZ[K_n]\otimes \ZZ[H_n]$ defined by
$$
   \Phi(b)= b\:(\si\circ \varphi)(b)^{-1} \otimes \varphi(b).
$$
Let us fix such a set-section.

 Recall that the {\em augmentation ideal} of a group $G$ is defined to
be the two-sided ideal $I(G)$ of $\ZZ[G]$ generated by the set 
$\{1-g; \; g\in G \}$. In this subsection, we prove the following.

\vspace{.3cm}
\begin{proposition}\label{proIJ} 
  The isomorphism $\Phi:\: \ZZ[B_n(M)]  \rightarrow \ZZ[K_n]\otimes 
\ZZ[H_n]$ sends isomorphically $J^d$ to $I(K_n)^d\otimes \ZZ[H_n]$
for all $d\in \NN$.
\end{proposition}
\vspace{.3cm}
  
Note that  Proposition~\ref{proIJ} implies that, in order to prove 
Theorem~\ref{prosepara}, it will suffice to prove the following two 
conditions.
\begin{enumerate}
 
 \item $\bigcap_{d=0}^{\infty} I(K_n)^d = \{ 0 \}$,

 \item $I(K_n)^d/I(K_n)^{d+1}\;$ is a free $\ZZ$-module for all $d\geq 1$.

\end{enumerate}

\vspace{.3cm}
 To prove Proposition~\ref{proIJ}, we will  make use of some classical 
exact sequences involving braid groups (see~\cite{birman}). 
The first one comes from the homomorphism $\pi$ 
which maps a given braid to the permutation that it induces on ${\cal P}$.
The kernel of this (clearly well defined) homomorphism is a subgroup 
of $B_n(M)$ denoted by $PB_n(M)$, whose elements are called {\em pure
braids}. Then one has:
$$
  1 \longrightarrow PB_n(M) \longrightarrow B_n(M) 
    \stackrel{\pi}{\longrightarrow} \Si_n \longrightarrow 1.
$$

On the other hand, there is a homomorphism $\ro: \; PB_n(M) \rightarrow 
PB_{n-1}(M)$ which sends  $(b_1,\ldots,b_n)$ to $(b_2,\ldots,b_n)$. 
If we set ${\cal P}_{n-1}=\{P_2,\ldots,P_n \}$, then the kernel of $\ro$ 
can be seen as the group $\pi_1(M\backslash {\cal P}_{n-1})$. This gives:
$$
  1 \longrightarrow \pi_1(M\backslash {\cal P}_{n-1}) 
    \longrightarrow PB_n(M) 
    \stackrel{\ro}{\longrightarrow}  PB_{n-1}(M)\longrightarrow 1.
$$

 Finally, if $b$ is a pure braid, the projection of each 
string $b_i$ ($i\in \{ 1,\ldots,n\}$) over $M$, denoted by 
$\overline{b_i}$, is a loop in $M$ based at $P_i$,
which determines an element $\mu_i\in \pi_1(M)$.
This gives a homomorphism $\theta:\:PB_n(M) \rightarrow\pi_1(M)^n$, which 
sends $(b_1,\ldots,b_n)$ to $(\mu_1,\ldots,\mu_n)$. 
One can easily verify that $K_n=\ker{\theta}$, and that the exact 
sequence
$$
  1 \longrightarrow K_n \longrightarrow PB_n(M) 
    \stackrel{\theta}{\longrightarrow} \pi_1(M)^n \longrightarrow 1
$$
extends to the exact sequence
$$
  1 \longrightarrow K_n \longrightarrow B_n(M) 
    \stackrel{\varphi}{\longrightarrow} H_n \longrightarrow 1.
$$
Moreover, $K_n$ is the normal closure in
$PB_n(M)$ of the subgroup $PB_n(D)$, where $D$ is a disc in 
$M$ which contains ${\cal P}$ (see \cite{birman}).

In what follows, we write $I=I(K_n)$ and we consider $\ZZ[K_n]$ as a
subring of $\ZZ[B_n(M)]$. The next lemma is a preliminary result
to the proof of Proposition~\ref{proIJ}.
 
\vspace{.3cm}
\begin{lemma}\label{bib}
 Let $B=\ZZ[B_n(M)]$. For every $d\geq 1$, one has
$$ 
    J^d = B \: I^d \: B = B I^d = I^d B.
$$
\end{lemma}
\vspace{.3cm}

\begin{proof}
  Since $K_n$ is a normal subgroup of $B_n(M)$, it is straightforward
to prove that $B \: I^d \: B = B I^d = I^d B$. So, it suffices to prove 
that $J=B\: I\: B$.

  The inclusion $J\subset B\: I\: B$ is obvious, once we notice that 
$\si_i^2\in K_n$ and that
$\si_i - \si_i^{-1} = \si_i^{-1} (\si_i^2 -1) \in B\: I\: B$. 
For the other inclusion, we must prove that for all $p\in K_n$ one has
$p-1\in J$. Suppose that $p=p_1\: p_2$, with $p_1,p_2 \in K_n$; then
$p-1=p_1(p_2-1)+(p_1-1)$, so it suffices to show it for a set of 
generators of $K_n$. As we said before, $K_n$ is the normal closure of
$PB_n(D)$ in $PB_n(M)$, so a set of generators of $K_n$ consists on 
the elements of the form $\al \: b \: \al^{-1}$, where $\al \in PB_n(M)$
and $b \in PB_n(D)$.

  Take an element $\al \: b \: \al^{-1}$ as above. One has:
$\al \: b \: \al^{-1}-1 = \al \: (b -1)\: \al^{-1}$, so we only have to
show that $b-1\in J$ for $b\in PB_n(D)$. It is known (\cite{stanford}, 
Lemma 1.2) that $b-1$ belongs to the ideal of $\ZZ[B_n(D)]$ generated
by $\{\si_i-\si_i^{-1}; \; i=1,\ldots,n-1\}$. But the extension of this 
ideal to $\ZZ[B_n(M)]\supset \ZZ[B_n(D)]$ is precisely $J$, so
$b-1\in J$.
\end{proof}

\vspace{.6cm}
\noindent{\sc Proof of Proposition~\ref{proIJ}}.
 First, we show that $\Phi(J^d)\subset I^d\otimes \ZZ[H_n]$ for all 
$d\geq 1$. By Lemma~\ref{bib}, we know that $J^d= I^d\: B$,
thus $J^d$ is generated as a $\ZZ$-module by the elements of the form 
$(k_1-1)\cdots (k_d-1)b$, where $b\in B_n(M)$ and  $k_i\in K_n$ 
for $i=1,\ldots,d$. Now, the image of such an element by $\Phi$ is 
$(k_1-1)\cdots (k_d-1)\:b' \otimes \varphi(b)$, where 
$b'=b\:(\si\circ \varphi)(b)^{-1}$, which clearly belongs to
$I^d\otimes \ZZ[H_n]$.

 The inclusion $I^d\otimes \ZZ[H_n] \subset \Phi(J^d)$ follows from
the facts that $I^d\otimes \ZZ[H_n]$ is generated as a $\ZZ$-module
by the elements of the form $(k_1-1)\cdots (k_d-1)k \otimes \be$, 
where $k_1,\ldots,k_d,k\in K_n$ and $\be\in H_n$, and that such an element 
is the image by $\Phi$ of 
$(k_1-1)\cdots (k_d-1)k\:\si(\be)\in I^d\: B =J^d$. 
$\fbox{}$

\subsection{The structure of $K_n$}\label{secalmdirpro}
 
The goal of this subsection is to prove the following.

\begin{proposition}\label{prosemidire}
 For $n\geq 2$, there exists a free group $F_n$ such that
$ K_n= F_n \rtimes K_{n-1}$.
Moreover, the action of $K_{n-1}$ on the abelianization of $F_n$ is trivial.
\end{proposition}

\vspace{.6cm}
\noindent {\bf Remarks.}
\begin{itemize}
 \item[(i)] The notation $F_n$ may lead to some confusion; indeed, here 
$F_n$ is not a free group of rank $n$. It is actually of infinite rank.

 \item[(ii)] A direct consequence of Proposition~\ref{prosemidire} is
that $K_n$ can be expressed as an iterated semi-direct product of
(infinitely generated) free groups
$$
   K_n= F_n \rtimes (F_{n-1} \rtimes(\cdots\rtimes (F_3 \rtimes F_2)\cdots )).
$$
\end{itemize}
 
   Recall the exact sequences defined in the previous subsection. Since
$K_n$ is a subgroup of $PB_n(M)$, we can consider the image by $\ro$ 
of $K_n$. By definition, it is equal to $K_{n-1}$. If we 
denote $F_n = \ker{\ro} \cap K_n$, we obtain the following commutative
diagram, where all rows and columns are exact:
$$
  \begin{array}{rcccccl}
  &     1      &      &     1     &      &    1     &   \\
  &  \uparrow  &      &  \uparrow &      &  \uparrow &   \\
 1\rightarrow & \pi_1(M,P_1) & \rightarrow & \pi_1(M)^n & \rightarrow & 
			\pi_1(M)^{n-1}  & \rightarrow 1 \\
 &  \uparrow & & \hspace{.3cm} \uparrow \mbox{\scriptsize $\theta$} &   &
   \hspace{.3cm} \uparrow \mbox{\scriptsize $\theta$} \\
 1\rightarrow  & \pi_1(M\backslash {\cal P}_{n-1}) & \rightarrow
   & PB_n(M) &\stackrel{\ro}{\rightarrow} & PB_{n-1}(M) & \rightarrow 1 \\
  &  \uparrow   & &  \uparrow  &  & \uparrow  &   \\
1 \rightarrow & F_n & \rightarrow & K_n & \stackrel{\ro}{\rightarrow} 
&  K_{n-1}  & \rightarrow 1  \\
  &   \uparrow & &  \uparrow  &  & \uparrow  &   \\
 &     1      &      &     1     &      &    1     & 
  \end{array}
$$

Notice that $F_n$ is a free group, since it is a subgroup of 
$\pi_1(M\backslash {\cal P}_{n-1})$, which is a free group. 
We are specially interested in the lowest row of the diagram. 
In particular, in order
to show Proposition~\ref{prosemidire}, we will
show that there exists a homomorphism  $s:\; K_{n-1} \rightarrow K_n$ 
which is a section of $\ro$, and that $K_{n-1}$ 
acts trivially on the abelianization of $F_n$.  
We turn first to find a free set of generators for $F_n$.

 Let $\Om=\{ \om_1,\ldots, \om_{2g}\}$ be a set of $2g$ letters.
It is well known that a presentation for $\pi_1(M)$ is as follows. 
$$
   \pi_1(M)=\left< \Om ;
\; \;(\om_1\om_2\cdots \om_{2g} \om_1^{-1}\om_2^{-1}\cdots 
\om_{2g}^{-1})=1\right>. 
$$

 For every element $\ga \in \pi_1(M)$ we
choose a unique word $\widetilde{\ga}$ over $\Om\cup \Om^{-1}$ 
which represents $\ga$. 
We call this word the {\em normal form of $\ga$}. 
Normal forms are chosen in such a way that they are prefix-closed (namely, 
if $\om_1\om_2$ is a normal form, then $\om_1$ is also a normal form). 
For every word $\om$ over $\Om\cup \Om^{-1}$, we will denote by 
$\om_{(i)}$ the word over $\{a_{i,1}^{\pm 1},\ldots,a_{i,2g}^{\pm 1}\}$ 
obtained from $\om$ 
by replacing $\om_j^{\pm 1}$ by $a_{i,j}^{\pm 1}$, for all $j=1,\ldots,2g$.

\begin{figure}[ht]
\centerline{ \begin{picture}(0,0)%
\includegraphics{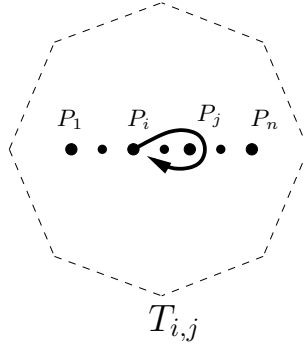}%
\end{picture}%
\setlength{\unitlength}{2960sp}%
\begingroup\makeatletter\ifx\SetFigFont\undefined
\def\x#1#2#3#4#5#6#7\relax{\def\x{#1#2#3#4#5#6}}%
\expandafter\x\fmtname xxxxxx\relax \def\y{splain}%
\ifx\x\y   
\gdef\SetFigFont#1#2#3{%
  \ifnum #1<17\tiny\else \ifnum #1<20\small\else
  \ifnum #1<24\normalsize\else \ifnum #1<29\large\else
  \ifnum #1<34\Large\else \ifnum #1<41\LARGE\else
     \huge\fi\fi\fi\fi\fi\fi
  \csname #3\endcsname}%
\else
\gdef\SetFigFont#1#2#3{\begingroup
  \count@#1\relax \ifnum 25<\count@\count@25\fi
  \def\x{\endgroup\@setsize\SetFigFont{#2pt}}%
  \expandafter\x
    \csname \romannumeral\the\count@ pt\expandafter\endcsname
    \csname @\romannumeral\the\count@ pt\endcsname
  \csname #3\endcsname}%
\fi
\fi\endgroup
\begin{picture}(2409,2703)(619,-2482)
\put(1801,-2401){\makebox(0,0)[lb]{\smash{\SetFigFont{14}{16.8}{rm}$T_{i,j}$}}}
\put(1036,-691){\makebox(0,0)[lb]{\smash{\SetFigFont{9}{10.8}{rm}$P_{1}$}}}
\put(2656,-691){\makebox(0,0)[lb]{\smash{\SetFigFont{9}{10.8}{rm}$P_{n}$}}}
\put(1621,-691){\makebox(0,0)[lb]{\smash{\SetFigFont{9}{10.8}{rm}$P_{i}$}}}
\put(2206,-646){\makebox(0,0)[lb]{\smash{\SetFigFont{9}{10.8}{rm}$P_{j}$}}}
\end{picture}
}
\caption{The path (or braid) $T_{i,j}$.}
\label{Tij}
\end {figure}

  Let us consider, for $1\leq i < j \leq n$, the braid $T_{i,j}$ drawn in
Figure~\ref{Tij}. All its strings are trivial  except the $i$-th one 
which goes around the points $P_{i+1},\ldots,P_j$ and turns back to $P_i$.

 Notice that in $\pi_1(M\backslash{\cal P}_{n-1})$, viewed as 
a subgroup of $PB_n(M)$, one has  
$$
T_{1,n}=a_{1,1}\cdots a_{1,2g} a_{1,1}^{-1}\cdots a_{1,2g}^{-1}.
$$

\vspace{.3cm}
\begin{lemma}\label{baseFn}
 The following set is a free system of generators for $F_n$.
$$
  {\cal B} = \{ \widetilde{\ga}_{(1)}\: T_{1,j}\: 
\widetilde{\ga}_{(1)}^{-1}\; ; \;\;  
2\leq j \leq n \mbox{ and }\ga\in \pi_1(M)\}.
$$

\end{lemma}

\vspace{.3cm}
\begin{proof}
Consider the Cayley graph of $\pi_1(M)$, which is defined as follows. 
Its vertices are the elements of $\pi_1(M)$, and its edges are 
labeled by $\Om$. For every vertex $\ga\in \pi_1(M)$ and for every
$i\in{1,\ldots,2g}$, there is exactly one edge labeled by $\om_i$, 
with source $\ga$ and target $\ga \om_i$.
 
 In this graph, the normal form of an element $\ga \in \pi_1(M)$
corresponds to a unique path going from $1$ to $\ga$. By the prefix-closed
condition mentioned above, the set of normal forms of $\pi_1(M)$
defines a maximal tree $T$ of the Cayley graph.

 The Cayley graph of $\pi_1(M)$ can be seen as the one-skeleton of 
a tiling of the plane. For every vertex $\ga$, the 
path which starts at $\ga$ and which is labeled by
$\om_1\ldots \om_{2g} \om_1^{-1}\ldots \om_{2g}^{-1}$ bounds a 
fundamental region $R_{\ga}$ of this tiling, and all fundamental regions
are obtained in this way. Hence, there is a one to one 
correspondence between the vertices of the Cayley graph and its fundamental
regions. Therefore, the fundamental group of the Cayley graph of $\pi_1(M)$
is the free group with free system of generators 
$$
  \{ \widetilde{\ga}\: (\om_1\ldots \om_{2g} \om_1^{-1}\ldots \om_{2g}^{-1})\: 
  \widetilde{\ga}^{-1}\; ; \; \;\; \ga\in \pi_1(M,P_1)\}.
$$

We now define a graph $\Ga$ as follows. Take the Cayley graph of 
$\pi_1(M)$ and replace the labels $\om_i$ by $a_{1,i}$.
Then, for every vertex $\ga$, add $n-2$ edges with
source and target $\ga$, labeled by $T_{1,2}, \ldots, T_{1,n-1}$, 
respectively. Notice that the fundamental group of $\Ga$ is the free
group with free system of generators ${\cal B}$, where 
$T_{1,n}=a_{1,1}\cdots a_{1,2g}a_{1,1}^{-1}\cdots a_{1,2g}^{-1}$.

 Recall the exact sequence
$$
  1 \longrightarrow F_n \longrightarrow \pi_1(M\backslash {\cal P}_{n-1})
  \stackrel{\theta}{\longrightarrow}  \pi_1(M,P_1) \longrightarrow 1.
$$
One can easily verify that 
$\pi_1(M\backslash {\cal P}_{n-1})$ is freely generated by the set
$$
  \left\{a_{1,1},\ldots, a_{1,2g}, T_{1,2},\ldots, T_{1,n-1}  \right\},
$$
and that $\theta$ sends $a_{1,i}$ to $\om_i$ for all $i=1,\ldots,2g$,
and sends $T_{1,j}$ to $1$ for all $j=2,\ldots,n-1$. It follows from 
classical geometric methods (see \cite{lyndonschupp}) that 
the group $F_n$ is the fundamental group of the graph $\Ga$, hence 
${\cal B}$ is a free system of generators for $F_n$, as we wanted to prove.
\end{proof}

\vspace{.6cm}
\begin{lemma}
 There is a homomorphism $\si:\: K_{n-1} \rightarrow K_n$ which is a section
of $\ro:\: K_n \rightarrow K_{n-1}$.
\end{lemma}

\begin{proof}
 The case $n=2$ is trivial, since 
$K_1=\ker{\left(\pi_1(M)\stackrel{\theta}{\rightarrow} \pi_1(M)\right)}=1$.
Hence $K_2=F_2$ is a free group of infinite rank.
Suppose now that $n>2$.

 It is well known that the kernel of the homomorphism 
$\theta_2:\:PB_n(M)\rightarrow \pi_1(M,P_2)$ 
is $PB_{n-1}(M\backslash \{P_2\})$ (see \cite{birman}). Moreover, one
can easily see that $K_n$ lies in this kernel, namely
$K_n\subset PB_{n-1}(M\backslash \{P_2\})$. Similarly, one has
$K_{n-1}\subset PB_{n-2}(M\backslash \{P_2\})$. The homomorphism
$\ro_n:\: PB_{n-1}(M\backslash \{P_2\})\rightarrow 
PB_{n-2}(M\backslash \{P_2\})$ which sends $(b_1,b_3,\ldots, b_n)$ 
to $(b_3,\ldots,b_n)$ is the restriction of $\ro$ to 
$PB_{n-1}(M\backslash \{P_2\})$. In particular, it sends $K_n$ onto
$K_{n-1}$. 

 We consider an embedding $f:\: M\backslash  \{P_2\} \rightarrow 
M\backslash \{P_2\}$ satisfying:
\begin{itemize}
 \item $f(P_i)=P_i$ for $i=3,\ldots,n$;

 \item $P_1$ does not lie in the image of $f$; 

 \item $f$ is a homotopy equivalence relatively to $\{P_3,\ldots,P_n\}$.
\end{itemize}

 Then, the map $f$ induces a homomorphism
$\si:\: PB_{n-2}(M\backslash \{P_2\}) \rightarrow 
PB_{n-1}(M\backslash \{P_2\})$, which sends $(b_3,\ldots,b_n)$ to
$(1_{P_1},(f\times \id)b_3, \ldots (f\times \id)b_n)$. By the third 
condition, this homomorphism is a section of $\ro_n$. It obviously sends
$K_{n-1}$ to $K_n$. 
\end{proof}

\vspace{.6cm}
  Now, $K_{n-1}$ acts on $F_n$ in the following way: Given $b\in K_{n-1}$,
the action induced by $b$ sends $f\in F_n$ to $\si(b) f \si(b)^{-1}$.
This action induces an action of $K_{n-1}$ on the abelianization 
$F_n/[F_n,F_n]$ of $F_n$ (here $[F_n,F_n]$ denotes the commutator subgroup 
of $F_n$). The proof of Proposition~\ref{prosemidire} is finally obtained 
from the following result.

\vspace{.3cm}
\begin{lemma}
 The action of $K_{n-1}$ on the abelianization of $F_n$ is trivial. 
\end{lemma}

\begin{proof}
  We only need to verify that the action of the generators of $K_{n-1}$
on the generators of $F_n$ is trivial after abelianization. 
Moreover, let us see that it 
suffices to show the result for the action defined by any set-map 
section $s$ of $\ro$.

  Indeed, if $s$ is a section of $\ro$, then for every
$b \in K_{n-1}$, there exists an element $\widehat{b}\in F_n$
such that $\si(b)=\widehat{b}\:s(b)$. Therefore, if $K_{n-1}$ acts
trivially on $F_n/[F_n,F_n]$ via $s$, we obtain, for every $f\in F_n$,
$$
  \si(b)\; f \;\si(b)^{-1} \equiv \widehat{b}\;\left( s(b)\; f \; 
s(b)^{-1}\right) 
\widehat{b}^{-1} \equiv \widehat{b} \; f \; \widehat{b}^{-1} \equiv 
f  \quad \quad \mbox{(mod $[F_n,F_n]$)}.
$$

  As we said in Subsection~\ref{secIJ}, a set of generators for $K_{n-1}$ 
consists on the elements of the form $\al \: b \: \al^{-1}$, where $\al 
\in PB_{n-1}(M)$ and $b \in PB_{n-1}(D)$. On the other hand, it is known 
that $T=\left\{T_{i,j}\; | \quad 2\leq i<j\leq n\right\}$ is a set of 
generators for $PB_{n-1}(D)$, where $T_{i,j}$ denotes the braid defined
in Subsection~\ref{coincides}. Therefore, the following is a set of 
generators for $K_{n-1}$:
$$
  \left\{ \al \: T_{i,j}\: \al^{-1}; \; \; 2\leq i<j \leq n, \; 
  \mbox{and }\al \mbox{ is a word over } \{a_{k,r}^{\pm 1}; \;2\leq k \leq n 
 \mbox{ and } 1\leq r \leq 2g \} \right\}.
$$
We take $s$ such that 
$s(\al \: T_{i,j}\: \al^{-1})=\al \: T_{i,j}\: \al^{-1}\in K_n$. 
In other words, we just add a trivial string based at $P_1$ for
any element of the set of generators of $K_{n-1}$. We remark that $s$ is
a set map section of $\ro$, but it is not a homomorphism.

  Now, $F_n$ is by definition a normal subgroup of $PB_n(M)$.  Therefore,
if we show that each $T_{i,j}$ $\; (i\geq 2)$ acts trivially on 
$F_n/[F_n,F_n]$ by conjugation, then we will have finished the proof, 
since in that case:
$$
  (\al \: T_{i,j}\: \al^{-1}) \: f \: (\al \: T_{i,j}^{-1}\: \al^{-1})
\equiv \al \: T_{i,j}\: (\al^{-1}  f \: \al) \: T_{i,j}^{-1}\: \al^{-1}
\equiv  \al \: (\al^{-1}   f \:\al) \: \al^{-1} \equiv f 
\quad \mbox{(mod $[F_n,F_n]$)}.
$$

  Recall from Lemma~\ref{baseFn} that 
$$
 {\cal B} = \{ \widetilde{\ga}_{(1)}\: T_{1,k}\: \widetilde{\ga}_{(1)}^{-1}\; ; \;\;  
2\leq k \leq n  \mbox{ and }\; \ga\in \pi_1(M)\}
$$
is a free system of generators for $F_n$, where $\widetilde{\ga}_{(1)}$ 
is a  word over $\{a_{1,1}^{\pm 1},\ldots,a_{1,2g}^{\pm 1}\}$ for all 
$\ga \in\pi_1(M) $. 
One can verify (just drawing the corresponding 
braids) that in $PB_n(M)$ one has the following relations:
$$
\begin{array}{lr}

  T_{i,j} a_{1,r} T_{i,j}^{-1} = a_{1,r}, & 
\\ \\
 T_{i,j}T_{1,k}T_{i,j}^{-1} = T_{1,k}  & 
(k<i \mbox{ or } k\geq j),
\\ \\

T_{i,j}T_{1,k}T_{i,j}^{-1}= T_{1,i-1}T_{1,i}^{-1}T_{1,k}T_{1,j}^{-1}
   T_{1,i}T_{1,i-1}^{-1}T_{1,j}\equiv T_{1,k} \quad 
(\mbox{mod }[F_n,F_n])  & (i\leq k < j),
\end{array}
$$
where $2\leq i < j \leq n$, $\; r\in \{1,\ldots, 2g \}$ and 
$k\in \{2,\ldots, n \}$.

 Therefore, in $F_n/[F_n,F_n]$, we have:
$$
   T_{i,j} \: (\widetilde{\ga}_{(1)}\: T_{1,k}\: 
\widetilde{\ga}_{(1)}^{-1})\: T_{i,j}^{-1} =
   \widetilde{\ga}_{(1)} \: (T_{i,j}\: T_{1,k}\: T_{i,j}^{-1}) 
\widetilde{\ga}_{(1)}^{-1} \equiv 
   \widetilde{\ga}_{(1)}\: T_{1,k}\: \widetilde{\ga}_{(1)}^{-1}
\quad \quad \mbox{(mod $[F_n,F_n]$)},
$$
as we wanted to show.
\end{proof}

\subsection{Proof of Theorem~\ref{prosepara}}\label{secsepara}

 Let $A$ and $C$ be two groups such that $C$ acts on $A$.
For $a\in A$ and $c\in C$, we denote by $a^c$ the action of $c$ on $a$.
Then, the $\ZZ$-module $\ZZ[A]\otimes \ZZ[C]$ carries a natural structure
of $\ZZ$-algebra, where the multiplication is defined by 
$$
  (a_1\otimes c_1)\cdot (a_2\otimes c_2) = (a_1\: a_2^{c_1})\otimes
 (c_1c_2). 
$$
Moreover, this algebra is naturally isomorphic to $\ZZ[A\rtimes C]$ via
an isomorphism which sends $a\otimes c$ to $ac$ for all $a\in A$ and 
all $c\in C$. 

 Recall that the augmentation ideal of a group $G$ is denoted by $I(G)$.
The following lemma will be used to prove Theorem~\ref{prosepara}.
Its proof can be found in [P, Lemma 3.1].

\vspace{.3cm}
\begin{lemma}\label{lempapa}
  Let $A$ and $C$ be two groups. Assume that an action of $C$ on $A$ is
given, and that this action induces the trivial action on the 
abelianization of $A$. Then one has
$$
  I(A\rtimes C)^m = \sum_{k=0}^m{I(A)^k \otimes I(C)^{m-k}} 
$$
for all $m\geq 0$. $\fbox{}$
\end{lemma}

\vspace{.6cm}

\noindent {\sc Proof of Theorem~\ref{prosepara}}.
As pointed out in Subsection~\ref{secIJ}, 
it suffices to prove the following two 
conditions.
\begin{enumerate}
 
 \item $\bigcap_{d=0}^{\infty} I(K_n)^d = \{ 0 \}$,

 \item $I(K_n)^d/I(K_n)^{d+1}\;$ is a free $\ZZ$-module for all $d\geq 0$.

\end{enumerate}

 We argue by induction on $n$. The case $n=1$ is
trivial, since $I(K_1)=I(\{1\})=0$. So, we assume that
$n\geq 2$ and that Conditions 1 and 2 hold for $K_p$, where $p<n$.

The group $F_n$ is free, thus, by \cite{fox}, one has
\begin{enumerate}
 
 \item $\bigcap_{d=0}^{\infty} I(F_n)^d = \{ 0 \}$,

 \item $I(F_n)^d/I(F_n)^{d+1}\;$ is a free $\ZZ$-module for all $d\geq 0$.

\end{enumerate}
These two properties imply that $I(F_n)^d$ is a free $\ZZ$-module and that 
$$
  I(F_n)^d\simeq \bigoplus_{k=d}^{\infty}{I(F_n)^k/I(F_n)^{k+1}}
$$
for all $d\geq 0$. We choose a $\ZZ$-basis ${\cal B}_d$ of
$I(F_n)^d/I(F_n)^{d+1}$ for all $d\geq 0$. From the above isomorphism one 
has that ${\cal B}_{\geq d}= \coprod_{k=d}^{\infty}{\cal B}_k$ is a
$\ZZ$-basis of $I(F_n)^d$. 

 From the induction hypothesis, one also has 
\begin{enumerate}
 
 \item $\bigcap_{d=0}^{\infty} I(K_{n-1})^d = \{ 0 \}$,

 \item $I(K_{n-1})^d/I(K_{n-1})^{d+1}\;$ is a free $\ZZ$-module for all 
$d\geq 0$.

\end{enumerate}
This implies that $I(K_{n-1})^d$ is a free $\ZZ$-module and that
$$
  I(K_{n-1})^d\simeq \bigoplus_{k=d}^{\infty}{I(K_{n-1})^k/I(K_{n-1})^{k+1}}
$$
for all $d\geq 0$. We choose a $\ZZ$-basis ${\cal C}_d$ of
$I(K_{n-1})^d/I(K_{n-1})^{d+1}$ for all $d\geq 0$. 
Thus ${\cal C}_{\geq d}= \coprod_{k=d}^{\infty}{\cal C}_k$ is a
$\ZZ$-basis of $I(K_{n-1})^d$. 

 Now, by Proposition~\ref{prosemidire}, $F_n$ and $K_{n-1}$ satisfy the
hypothesis of Lemma~\ref{lempapa}. Hence, we have the equality 
$$
  I(K_n)^m= \sum_{d=0}^m{I(F_n)^d \otimes I(K_{n-1})^{m-d}} 
$$
for all $m\geq 0$. From this equality, one can easily verify that the set
$$
   {\cal D}_{\geq m} = \{ b\otimes c \in \ZZ[F_n]\otimes\ZZ[K_{n-1}]
 \; ; \quad b\in {\cal B}_i, \; c\in {\cal C}_j, \; i+j\geq m \}
$$
is a generating set for $I(K_n)^m$. Since this set is linearly 
independent, $I(K_n)^m$ is a free $\ZZ$-module whose basis is 
${\cal D}_{\geq m}$. It follows that 
$I(K_n)^m/I(K_n)^{m+1}$ is a free $\ZZ$-module with basis
$$
   {\cal D}_{m} ={\cal D}_{\geq m}\backslash{\cal D}_{\geq m+1} =
\{ b\otimes c \in \ZZ[F_n]\otimes\ZZ[K_{n-1}]
 \; ; \quad b\in {\cal B}_i, \; c\in {\cal C}_j, \; i+j= m \},
$$
and that $\bigcap_{d=0}^{\infty}I(K_n)^d=\{0\}$, since 
${\cal D}_{\geq 0}$ is a basis for $\ZZ[K_n]$.
$\fbox{}$

\section{The universal Vassiliev invariant}

 The proof of Theorem~\ref{teorisom} is divided into five steps. In what 
follows each subsection will correspond to one of them.

 The first subsection is dedicated to the definition of a linear map
$$
 u: \: \ZZ[B_n(M)] \longrightarrow \widehat{\cal A}_n \rtimes \ZZ[H_n].
$$
Recall that $\ZZ[B_n(M)]$ is isomorphic to $\ZZ[K_n]\otimes \ZZ[H_n]$ as 
a $\ZZ$-module (see Subsection~\ref{secIJ}), and notice that
$\widehat{\cal A}_n \rtimes \ZZ[H_n]$ is equal 
as a $\ZZ$-module to $\widehat{\cal A}_n \otimes \ZZ[H_n]$. 
Hence, we will only need to define a linear map $v:\ZZ[K_n]\rightarrow 
\widehat{\cal A}_n$.

 Recall that the subgroups $G_i$ of the {\it lower central series} of
a group $G$ are defined recursively by 
$G_1=G$ and $G_{i+1}= [G,G_i]$ for $i\geq 1$. 
The {\em associated graded Lie algebra of $G$}
is defined by $\mbox{gr}(G) = \bigoplus_{i\geq1}{G_i/G_{i+1}}$. 
It is a graded Lie algebra over $\ZZ$
whose enveloping algebra is denoted by ${\cal U}{\mbox{gr}(G)}$. In 
subsection~\ref{chi1} we construct a homomorphism 
$\chi_1:\: {\cal A}_n \rightarrow {\cal U}{\mbox{gr}(K_n)}$ of
$\ZZ$-algebras and we prove that this homomorphism is actually an
isomorphism.

 Write $I=I(K_n)$. In Subsection~\ref{chi2} we construct an isomorphism
$\chi_2:\: {\cal U}{\mbox{gr}(K_n)}\rightarrow \mbox{gr}_I\ZZ[K_n]$
of $\ZZ$-algebras using some Quillen's results.  
 
 In Subsection~\ref{inverso} we consider the isomorphism 
$\chi=\chi_2\circ \chi_1: \: {\cal A}_n\rightarrow \mbox{gr}_I\ZZ[K_n]$
and we prove that gr$v$ is the inverse of $\chi$. It will immediately
follow that 
$\mbox{gr}v: \: \mbox{gr}_I\ZZ[K_n]\rightarrow {\cal A}_n$ is an 
isomorphism of $\ZZ$-algebras, and that 
$\mbox{gr}u: \: \mbox{gr}_V\ZZ[B_n(M)]\rightarrow {\cal A}_n\rtimes \ZZ[H_n]$
is an isomorphism of $\ZZ$-modules.

 Finally, we prove in Subsection~\ref{gruhomo} that gr$u$ is a homomorphism 
of $\ZZ$-algebras. This will finish the proof of 
Theorem~\ref{teorisom}.

\subsection{Construction of $u: \: \ZZ[B_n(M)] \longrightarrow 
\widehat{\cal A}_n \rtimes \ZZ[H_n]$}

 From now on we fix a set-section  $\si:\: H_n\rightarrow B_n(M)$ of 
$\varphi:\:  B_n(M)\rightarrow H_n$. As we pointed out in 
Subsection~\ref{secIJ}, this set-section leads to an isomorphism
$\Phi:\: \ZZ[B_n(M)] \longrightarrow \ZZ[K_n]\otimes \ZZ[H_n]$ 
of $\ZZ$-modules.

 We turn now to define a linear homomorphism
$v:\ZZ[K_n] \rightarrow \widehat{\cal A}_n$. Then we will set
$$
  u= (v\otimes \mbox{id})\circ \Phi: \: 
  \ZZ[B_n(M)] \longrightarrow \widehat{\cal A}_n \rtimes \ZZ[H_n]
  =\widehat{\cal A}_n \otimes \ZZ[H_n].
$$

\begin{figure}[ht]
\centerline{\begin{picture}(0,0)%
\includegraphics{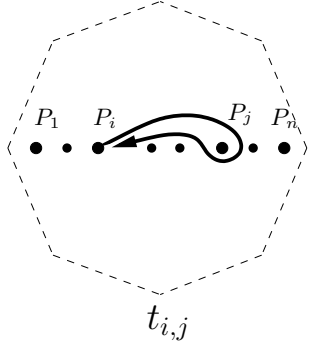}%
\end{picture}%
\setlength{\unitlength}{2960sp}%
\begingroup\makeatletter\ifx\SetFigFont\undefined
\def\x#1#2#3#4#5#6#7\relax{\def\x{#1#2#3#4#5#6}}%
\expandafter\x\fmtname xxxxxx\relax \def\y{splain}%
\ifx\x\y   
\gdef\SetFigFont#1#2#3{%
  \ifnum #1<17\tiny\else \ifnum #1<20\small\else
  \ifnum #1<24\normalsize\else \ifnum #1<29\large\else
  \ifnum #1<34\Large\else \ifnum #1<41\LARGE\else
     \huge\fi\fi\fi\fi\fi\fi
  \csname #3\endcsname}%
\else
\gdef\SetFigFont#1#2#3{\begingroup
  \count@#1\relax \ifnum 25<\count@\count@25\fi
  \def\x{\endgroup\@setsize\SetFigFont{#2pt}}%
  \expandafter\x
    \csname \romannumeral\the\count@ pt\expandafter\endcsname
    \csname @\romannumeral\the\count@ pt\endcsname
  \csname #3\endcsname}%
\fi
\fi\endgroup
\begin{picture}(2409,2703)(619,-2482)
\put(2836,-691){\makebox(0,0)[lb]{\smash{\SetFigFont{9}{10.8}{rm}$P_{n}$}}}
\put(1801,-2401){\makebox(0,0)[lb]{\smash{\SetFigFont{14}{16.8}{rm}$t_{i,j}$}}}
\put(1351,-691){\makebox(0,0)[lb]{\smash{\SetFigFont{9}{10.8}{rm}$P_{i}$}}}
\put(2476,-646){\makebox(0,0)[lb]{\smash{\SetFigFont{9}{10.8}{rm}$P_{j}$}}}
\put(856,-691){\makebox(0,0)[lb]{\smash{\SetFigFont{9}{10.8}{rm}$P_{1}$}}}
\end{picture}
}
\caption{The braid $t_{i,j}$.}
\label{tij}
\end {figure}

Recall that, for all $\ga\in \pi_1(M)$ and all $i\in \{ 1,\ldots,n\}$, 
we denote by $\widetilde{\ga}_{(i)}$ the normal form of $\ga$ over the 
generators $\{ a_{i,1}^{\pm 1}, \ldots, a_{i,2g}^{\pm 1}\}$ of 
$\pi_1(M,P_i)$. For $1\leq i<j \leq n$, we write 
$t_{i,j}=t_{j,i}= T_{i,j}T_{i,j-1}^{-1}$,
which is the braid drawn in Figure~\ref{tij}. These braids are the classical 
generators of $PB_n(D)$. Then, for $i\neq j$,
we denote by $f_{i,j,\ga}$ the element $\widetilde{\ga}_{(i)} t_{i,j} 
\widetilde{\ga}_{(i)}^{-1}$ of $B_n(M)$. From Lemma~\ref{baseFn} 
follows that  $F_{(n+1)-i}$ is the free group freely generated by
$$
   {\cal F}_{i,n} = \{f_{i,j,\ga}\; ; \;\;  
         j=i+1,\ldots,n, \; \; \ga\in \pi_1(M)\}.
$$
 Moreover, it is shown in Subsection~\ref{secalmdirpro} that 
$K_n=(F_n \rtimes (F_{n-1} \rtimes (\cdots \rtimes(F_3 \rtimes F_2)\cdots))$.
So, every element $k\in K_n$ can be uniquely written in the form 
$k=k_1\cdots k_{n-1}$, where $k_i$ is a reduced word 
over ${\cal F}_{i,n}\cup {\cal F}_{i,n}^{-1}$.
Now, for $i\in \{1,\ldots,n-1 \}$, there is an injective multiplicative
homomorphism $u_i:\: F_{(n+1)-i} \rightarrow \ZZ\lbra t_{i,j,\ga}\rbra$,
where $\ZZ\lbra t_{i,j,\ga}\rbra$ denotes the ring of non commutative
formal power series over non-commutative variables $t_{i,j,\ga}$, defined by 
$$
\begin{array}{l}
    u_i(  f_{i,j,\ga})= 1+ t_{i,j,\ga} \\
    u_i( f_{i,j,\ga}^{-1}) = 1 - t_{i,j,\ga} + t_{i,j,\ga}^2 - \cdots
\end{array}
$$
This well known homomorphism is called {\em Magnus expansion} of 
$F_{(n+1)-i}$ (see \cite{magnus}).
We denote by  $v_i$ the composition of $u_i$ with the
canonical projection $\ZZ \lbra t_{i,j,\ga}\rbra \rightarrow 
\widehat{\cal A}_n$, and we finally define the linear map
$ v:\: \ZZ[K_n]  \rightarrow  \widehat{\cal A}_n$ by
$$
 v(k) = v_1(k_1) v_2(k_2) \cdots v_{n-1}(k_{n-1}),
$$
where $k=k_1 k_2\cdots k_{n-1}$ is the decomposition of $k \in K_n$
defined above.


\subsection{The isomorphism $\chi_1:\: {\cal A}_n \rightarrow 
{\cal U}{\mbox{gr}(K_n)}$}\label{chi1}

The goal of this subsection is to prove the following.

\begin{proposition}\label{prochi1}
There is a well defined isomorphism $\chi_1:\: {\cal A}_n \rightarrow 
{\cal U}{\mbox{\em gr}(K_n)}$ of $\ZZ$-algebras which sends $t_{i,j,\ga}$
to $f_{i,j,\ga}$ for all $i,j\in \{1,\ldots,n \}$, $i\neq j$, and all
$\ga\in \pi_1(M)$.
\end{proposition}

\vspace{.3cm}
 Consider the graded Lie algebra $L_n$ given by the following presentation:
\begin{itemize}
  \item {\bf Generators:} \quad $\{t_{i,j,\ga}; \; \; 1\leq i,j \leq n,
		\; i\neq j, \; \ga\in \pi_1(M) \}$.

  \item {\bf Relations:}
	\begin{itemize}
	  
	  \item[(L1)] $t_{i,j,\ga}= t_{j,i,\ga^{-1}}$, \quad
		for all $i,j\in\{1,\ldots,n\}$, $i\neq j$, and all 
		$\ga\in \pi_1(M)$,

	  \item[(L2)] $[t_{i,j,\ga}\: ,\: t_{k,l,\de}]=0,$ \quad 
  for all distinct  $i,j,k,l\in\{1,\ldots, n\}$ and all 
  $\ga, \de \in \pi_1(M)$,

          \item[(L3)] $[t_{i,j,\ga}\: ,\: t_{j,k,\de}+t_{i,k,(\ga\de)}]=0,\;$ 
	        for all distinct $i,j,k \in\{1,\ldots, n\}$ 
		and all $\ga, \de \in \pi_1(M)$,

	\end{itemize}
\end{itemize}
where $[\_\:,\_]$ denotes the Lie bracket. 
 
 One has ${\cal U}{L_n}={\cal A}_n$, so, in order to prove 
Proposition~\ref{prochi1}, it suffices to prove the following.

\vspace{.3cm}
\begin{proposition}\label{isolie}
There is a well defined Lie algebra isomorphism  
$\psi_n:\: L_n \longrightarrow  \mbox{\em gr}(K_n)$ 
which sends $t_{i,j,\ga}$ to $f_{i,j,\ga}$ for all $i,j\in \{1,\ldots,n \}$, 
$i\neq j$, and all $\ga\in \pi_1(M)$.
\end{proposition}

\vspace{.3cm}
 The following lemmas \ref{lemar1} to \ref{biendef} are preliminary 
results to the proof of Proposition~\ref{isolie}.

\begin{lemma}\label{lemar1}
  Let $\om$ be a word over $\Om^{\pm 1}$. Then there exists 
$W_{\om}\in (K_n)_2=[K_n,K_n]$ such that
$$
  \om_{(j)}\: t_{i,j} \:\om_{(j)}^{-1}= (\om_{(i)}^{-1} \:t_{i,j}\: 
  \om_{(i)})\: W_{\om}.
$$
\end{lemma}

\begin{proof}
  We can suppose, without loss of generality, that $i<j$. Suppose first 
that $\om$ is a single letter. If $\om_{(j)}= a_{j,r}$ and
$r$ is odd, then one can easily show by drawing the braids that the 
following equality holds in $PB_n(M)$:
$$
  a_{j,r}\: t_{i,j}\: a_{j,r}^{-1} = (t_{i,j-1} \cdots t_{i,i+1})
 \: a_{i,r}^{-1} \: t_{i,j} \: a_{i,r}\: 
 (t_{i,i+1}^{-1} \cdots t_{i,j-1}^{-1}). 
$$
Hence
$$
a_{j,r}\: t_{i,j}\: a_{j,r}^{-1}=(a_{i,r}^{-1} \:t_{i,j}\: a_{i,r})
\: W_{\om},
$$
where
$$
W_{\om}= \left[a_{i,r}^{-1}\: t_{i,j}^{-1}\: 
a_{i,r} \; , \; t_{i,j-1} \cdots t_{i,i+1} \right] \in (K_n)_2.   
$$
 
  If $\om_{(j)}= a_{j,r}$ and $r$ is even, then one has 
$$
  a_{j,r}\: t_{i,j}\: a_{j,r}^{-1}  = a_{i,r}^{-1} \: (t_{i,i+1}^{-1} \cdots 
t_{i,j-1}^{-1} \:t_{i,j} \cdots t_{i,i+1})\: a_{i,r}.
$$
Therefore, 
$$
a_{j,r}\: t_{i,j}\: a_{j,r}^{-1} =(a_{i,r}^{-1} \:t_{i,j}\: a_{i,r})
\: W_{\om},
$$
where
$$
W_{\om}= \left[a_{i,r}^{-1}\: t_{i,j}^{-1}\: 
a_{i,r} \; , \; a_{i,r}^{-1} \: (t_{i,i+1}^{-1} \cdots 
t_{i,j-1}^{-1}) \:a_{i,r} \right] \in (K_n)_2. 
$$

  The computations for $\om_{(j)}=a_{j,r}^{-1}$ are the same  as for
$\om_{(j)}=a_{j,r}$ interchanging the case $r$ odd with the case $r$ even.
 
  Suppose now that $\om$ is a word of length $k$, $k>1$, and that the
result is true for words of length less than $k$. We write
$\om=\al \: \be$, with $|\al| , |\be| < k$. Consider
$$
\begin{array}{l}
 W' =  \al_{(j)}^{-1}\: W_{\be}  \: \al_{(j)}, \\
 W'' = \left[ \be_{(i)}^{-1} \: \al_{(j)} \:
 t_{i,j}^{-1}\:\al_{(j)}^{-1}\:\be_{(i)}\: , \: [\al_{(j)}\:, \: 
 \be_{(i)}^{-1}]  \right] \: W', \\
 W = \be_{(i)}^{-1}\: W_{\al} \: \be_{(i)}\: W''. 
\end{array}  
$$

 The hypothesis $i\neq j$ implies that 
$[\al_{(j)}\:,\:\be_{(i)}^{-1}]\in K_n$. Furthermore, both $K_n$ and
$(K_n)_2$ are normal subgroups of $PB_n(M)$, thus $W\in (K_n)_2$.
Finally, a direct calculation shows that:
$$
  \om_{(j)}\: t_{i,j}\: \om_{(j)}^{-1}  = 
(\om_{(i)}^{-1} \: t_{i,j}\: \om_{(i)}) \: W,
$$
as we wanted to show.
\end{proof}

\vspace{.6cm}
 
\begin{lemma}\label{lemaprev1}
  Let $i,j\in \{2,\ldots,n\}$, $i\neq j$, let
$\ga\in \pi_1(M)$, and let $\om$ be a word over $\Om^{\pm 1}$. 
Then there exists $W\in (K_n)_2$ such that
$$
   f_{i,j,\ga}\: \om_{(1)} \: f_{i,j,\ga}^{-1} = W\: \om_{(1)}.
$$

\end{lemma}

\begin{proof}
  Recall the epimorphism $\ro:\: PB_n(M) \longrightarrow PB_{n-1}(M)$. 
One has $\ro(\widetilde{\ga}_{(i)}^{-1}\:\om_{(1)} \:
\widetilde{\ga}_{(i)} )=1$, thus $\widetilde{\ga}_{(i)}^{-1}\:\om_{(1)} 
\:\widetilde{\ga}_{(i)}= b$, where $b$ is a word
over ${\cal B}=\{a_{1,1}^{\pm 1},\ldots, a_{1,2g}^{\pm 1}, 
t_{1,2}^{\pm 1},\ldots, t_{1,n}^{\pm 1}\}$.

  By drawing the braids, one sees that 
$t_{i,j}\: a_{1,r}^{\pm 1} \:t_{i,j}^{-1}= a_{1,r}^{\pm 1}$, 
for $r=1,\ldots, 2g$. 
Moreover, since $t_{i,j}, t_{1,k}\in K_n$, one has  
$t_{i,j}\: t_{1,k}^{\pm 1} \: t_{i,j}^{-1} = W\: t_{1,k}^{\pm 1}$, where 
$W\in (K_n)_2$, for all $k=2,\ldots,n$. Therefore, since $b$ is a word 
over ${\cal B}$, one has: $t_{i,j}\: b \: t_{i,j}^{-1}= W_b \: b$, 
where $W_b\in(K_n)_2$. Hence,
\begin{eqnarray*}
  f_{i,j,\ga}\: \om_{(1)} \: f_{i,j,\ga}^{-1} & = & 
   \left( \widetilde{\ga}_{(i)}\: t_{i,j} \: \widetilde{\ga}_{(i)}^{-1} 
\right) \: \om_{(1)} \:
 \left( \widetilde{\ga}_{(i)}\: t_{i,j}^{-1} \: \widetilde{\ga}_{(i)}^{-1} \right)  \\
  & = &  \widetilde{\ga}_{(i)}\: t_{i,j} \: b \: t_{i,j}^{-1} \: \widetilde{\ga}_{(i)}^{-1} \\
 & = & \widetilde{\ga}_{(i)}\: W_b \: b \: \widetilde{\ga}_{(i)}^{-1} \\
 & = & W\:  \widetilde{\ga}_{(i)}\: b \: \widetilde{\ga}_{(i)}^{-1} \\
 & = & W \: \om_{(1)},
\end{eqnarray*}
where $W= \widetilde{\ga}_{(i)}\: W_b \: \widetilde{\ga}_{(i)}^{-1}
\in (K_n)_2$, as we wanted to show.
\end{proof}

\vspace{.6cm}
\begin{lemma}\label{lemaprev2}
  Let $\om$ be a word over $\Om^{\pm 1}$, and let $i,j,k\in \{1,\ldots,n\}$, 
all distinct. Then there exists $W\in (K_n)_2$ such that
$$
   \om_{(i)} \: t_{j,k} \: \om_{(i)}^{-1} = W\: t_{j,k}.
$$
\end{lemma}

\begin{proof}
  Clearly, it suffices to show the lemma when $\om$ is a single letter. 
Besides, we can suppose that $j<k$. 
Then the result is a consequence of the following relations in 
$PB_n(M)$.
$$
\begin{array}{lll}
  a_{i,r}\: t_{j,k} \: a_{i,r}^{-1} = t_{j,k} & 
a_{i,r}^{-1}\: t_{j,k}\: a_{i,r}= t_{j,k},
 &  \mbox{if $i<j$ or $i>k$.} \\  
 & &  \\
a_{i,r}\: t_{j,k} \: a_{i,r}^{-1} = t_{j,i}^{-1} \: t_{j,k} \: t_{j,i},
\hspace{.5cm}
& a_{i,r}^{-1}\: t_{j,k}\: a_{i,r}= \al   \: t_{j,k} \: \al^{-1},
&  \mbox{if $j<i<k$ and $r$ is odd.} \\
&  & \\
 a_{i,r}\: t_{j,k} \: a_{i,r}^{-1}= \al   \: t_{j,k} \: \al^{-1},  &
a_{i,r}^{-1}\: t_{j,k}\: a_{i,r}=t_{j,i}^{-1} \: t_{j,k} \: t_{j,i}, 
\hspace{.5cm} &
   \mbox{if $j<i<k$ and $r$ is even,}  
\end{array}
$$
where $\al= \left( a_{j,r}^{-1} \: t_{j,j+1}^{-1}\cdots t_{j,i-1}^{-1}\: 
t_{j,i}\cdots t_{j,j+1} \: a_{j,r}\right)\in K_n$.
These relations can be easily verified by drawing pictures.
\end{proof}

\vspace{.6cm}

\begin{lemma}\label{tt}
 Let $i,j,k,l\in\{1,\ldots,n \}$, all distinct. Then
$$
\begin{array}{ll}
  [t_{i,j}, t_{k,l}] \equiv 0 & (\mbox{mod }(K_n)_3), \\
 & \\   
 \hspace{0cm} [t_{i,j}, t_{i,k}] \equiv [t_{j,k}, t_{i,j}] 
\hspace{1cm} & (\mbox{mod }(K_n)_3). 
\end{array}
$$
\end{lemma}

\begin{proof}
 This lemma follows from the well known congruences
$$
\begin{array}{ll}
  [t_{i,j}, t_{k,l}] \equiv 0 & (\mbox{mod }(PB_n(D))_3), \\
 & \\   
 \hspace{0cm} [t_{i,j}, t_{i,k}] \equiv [t_{j,k}, t_{i,j}] 
\hspace{1cm} & (\mbox{mod }(PB_n(D))_3), 
\end{array}
$$
(see, for example,  \cite{bar-natan3}), together with the inclusion
$(PB_n(D))_3 \subset (K_n)_3$.
\end{proof}

\vspace{.6cm}
 
\begin{lemma}\label{biendef}
There is a well defined Lie algebra homomorphism  
$\psi_n:\: L_n \longrightarrow  \mbox{\em gr}(K_n)$ 
which sends $t_{i,j,\ga}$ to $f_{i,j,\ga}$ for all $i,j\in \{1,\ldots,n \}$, 
$i\neq j$, and all $\ga\in \pi_1(M)$.
\end{lemma}

\vspace{.3cm}
\noindent {\sc Proof}:  
We have to show that the following congruences hold:
\begin{tabbing}
     (R1)$\;$ \= $f_{i,j,\ga}\equiv f_{j,i,\ga^{-1}}$ \hspace{3cm} \=
                 $(\mbox{mod }(K_n)_2)$, \hspace{.3cm}  \=
		for all $i,j\in\{1,\ldots,n\}$, $i\neq j$,  \\
        \> \>  \> and all $\ga\in \pi_1(M)$;  \\
    (R2) \> $[f_{i,j,\ga}\: ,\: f_{k,l,\de}]\equiv 0 $ \>
       $(\mbox{mod }(K_n)_3)$,
  \> for all distinct  $i,j,k,l\in\{1,\ldots, n\}$  \\
  \> \> \> and all $\ga, \de \in \pi_1(M)$;  \\
    (R3) \> $[f_{i,j,\ga}\: ,\: f_{j,k,\de}]\equiv
                 [f_{i,k,(\ga\de)}\: , \: f_{i,j,\ga}]\;$
           \>      $(\mbox{mod }(K_n)_3)$, 
	   \>     for all distinct $i,j,k \in\{1,\ldots, n\}$ \\
	\> \> \> and all $\ga, \de \in \pi_1(M)$.
\end{tabbing}

\vspace{.3cm}
 Notice that (R1) follows from Lemma~\ref{lemar1}. So, it remains to prove
(R2) and (R3). We argue by induction on $n$. The conditions (R2) and (R3)
being empty if $n=2$, we may assume that $n>2$, that (R2) holds if
$i,j,k,l\in\{2, \ldots, n \}$ (by induction), and that (R3) holds if
$i,j,k \in\{2, \ldots, n \}$ (by induction).

 We turn now to prove (R2) for $k=1$. By Lemma~\ref{lemaprev1}, there
exists $W_1\in (K_n)_2$ such that 
$$
  f_{i,j,\ga} \: \widetilde{\de}_{(1)} \:  f_{i,j,\ga}^{-1} 
  = W_1 \: \widetilde{\de}_{(1)}.
$$
Also, by Lemma~\ref{lemaprev2}, there exists $W_2\in (K_n)_2$ such that 
$$
\widetilde{\ga}_{(i)}^{-1}\: t_{1,l} \: \widetilde{\ga}_{(i)} =
W_2 \: t_{1,l}.
$$
Then
\begin{eqnarray*}
  f_{i,j,\ga}\: f_{1,l,\de}\: f_{i,j,\ga}^{-1} & = & 
  \left(f_{i,j,\ga} \: \widetilde{\de}_{(1)} \:  f_{i,j,\ga}^{-1} \right) 
  \: f_{i,j,\ga} \: t_{1,l} \:  f_{i,j,\ga}^{-1} \: 
  \left(f_{i,j,\ga} \: \widetilde{\de}_{(1)}^{-1} \:  f_{i,j,\ga}^{-1} \right)  \\
  & = & W_1 \: \left( \widetilde{\de}_{(1)} \: f_{i,j,\ga} \: t_{1,l} \:  f_{i,j,\ga}^{-1}
     \: \widetilde{\de}_{(1)}^{-1}  \right) \: W_1^{-1} \\
 & \equiv &
 \widetilde{\de}_{(1)} \: f_{i,j,\ga} \: t_{1,l} \:  f_{i,j,\ga}^{-1} \: \widetilde{\de}_{(1)}^{-1}
\hspace{1cm} \mbox{(mod $(K_n)_3$)} \\
 & = & \widetilde{\de}_{(1)} \:\left(\widetilde{\ga}_{(i)} \: t_{i,j} 
\: \widetilde{\ga}_{(i)}^{-1}\right) \: t_{1,l} \:\left(\widetilde{\ga}_{(i)} \: 
t_{i,j}^{-1} \: \widetilde{\ga}_{(i)}^{-1}\right) \: \widetilde{\de}_{(1)}^{-1} \\
& =  & \widetilde{\de}_{(1)} \: \widetilde{\ga}_{(i)} \: t_{i,j} 
\: W_2 \: t_{1,l} \:  t_{i,j}^{-1} \: 
\widetilde{\ga}_{(i)}^{-1} \: \widetilde{\de}_{(1)}^{-1} \\
& \equiv & \widetilde{\de}_{(1)} \: \widetilde{\ga}_{(i)} 
\: W_2 \: t_{i,j} \: t_{1,l} \:  t_{i,j}^{-1} \: 
\widetilde{\ga}_{(i)}^{-1} \: \widetilde{\de}_{(1)}^{-1}  
\hspace{1cm} \mbox{(mod $(K_n)_3$)}\\
& \equiv & \widetilde{\de}_{(1)} \: \widetilde{\ga}_{(i)}  
\: W_2 \: t_{1,l} \: \widetilde{\ga}_{(i)}^{-1} \: 
\widetilde{\de}_{(1)}^{-1}  \hspace{1cm} \mbox{(mod $(K_n)_3$)}
\hspace{1cm} \mbox{(by Lemma~\ref{tt})} \\
& = & 
\widetilde{\de}_{(1)} \: t_{1,l} \: \widetilde{\de}_{(1)}^{-1} = f_{1,l,\de}.  
\end{eqnarray*}
Therefore, (R2) holds for $k=1$.

 The congruence (R2) holds for either $i=1$, or $j=1$, or $l=1$, because
of the above case and of the relation (R1).

 We turn now to prove (R3) for $i=1$. By Lemma~\ref{lemaprev1},
there exists $W_1\in (K_n)_2$ such that 
$$
  f_{j,k,\de} \: \widetilde{\ga}_{(1)} \:  f_{j,k,\de}^{-1} 
  = W_1 \: \widetilde{\ga}_{(1)}.
$$
Also, by Lemma~\ref{lemar1}, there exists $W_2\in (K_n)_2$ such that 
$$
\widetilde{\de}_{(j)}^{-1}\: t_{1,j} \: \widetilde{\de}_{(j)} =
\widetilde{\de}_{(1)} \: t_{1,j}\: \widetilde{\de}_{(1)}^{-1}\:  W_2.
$$
Then
\begin{eqnarray*}
  f_{j,k,\de} \: f_{1,j,\ga} \: f_{j,k,\de}^{-1} & = &
\left(f_{j,k,\de} \:\widetilde{\ga}_{(1)}  \: f_{j,k,\de}^{-1} \right)
\left(f_{j,k,\de} \: t_{1,j} \: f_{j,k,\de}^{-1} \right)
\left(f_{j,k,\de} \: \widetilde{\ga}_{(1)}^{-1}\: f_{j,k,\de}^{-1} \right)
\\
& = & W_1 \: \left( \widetilde{\ga}_{(1)} \: f_{j,k,\de} \: t_{1,j} 
\: f_{j,k,\de}^{-1} \: \widetilde{\ga}_{(1)}^{-1} \right) \: 
W_1^{-1} \\
& \equiv &  \widetilde{\ga}_{(1)} \: f_{j,k,\de} \: t_{1,j} 
\: f_{j,k,\de}^{-1} \: \widetilde{\ga}_{(1)}^{-1}
 \hspace{1cm} \mbox{(mod $(K_n)_3$)} \\
& = & \widetilde{\ga}_{(1)} \:\left(\widetilde{\de}_{(j)} \: t_{j,k} 
\: \widetilde{\de}_{(j)}^{-1}\right) \: t_{1,j} \:\left(\widetilde{\de}_{(j)} \: 
t_{j,k}^{-1} \: \widetilde{\de}_{(j)}^{-1}\right) \: \widetilde{\ga}_{(1)}^{-1} \\
& = & \widetilde{\ga}_{(1)} \:\widetilde{\de}_{(j)} \: t_{j,k} 
\: \widetilde{\de}_{(1)} \: t_{1,j}\: \widetilde{\de}_{(1)}^{-1}\:  W_2\: 
t_{j,k}^{-1} \: \widetilde{\de}_{(j)}^{-1} \: \widetilde{\ga}_{(1)}^{-1} \\
& \equiv & \widetilde{\ga}_{(1)} \:\widetilde{\de}_{(j)} \: \left( t_{j,k} 
\: \widetilde{\de}_{(1)} \: t_{1,j}\: \widetilde{\de}_{(1)}^{-1}\: 
t_{j,k}^{-1}\right) \: W_2\:  \widetilde{\de}_{(j)}^{-1} \: 
\widetilde{\ga}_{(1)}^{-1}  \hspace{1cm} \mbox{(mod $(K_n)_3$)}  \\
& = & \widetilde{\ga}_{(1)} \:\widetilde{\de}_{(j)} \: \left[ t_{j,k} 
\: , \:  \widetilde{\de}_{(1)} \: t_{1,j}\: \widetilde{\de}_{(1)}^{-1} 
\right] \left(\widetilde{\de}_{(1)} \: t_{1,j}\: \widetilde{\de}_{(1)}^{-1} 
\right) \: W_2\:  \widetilde{\de}_{(j)}^{-1} \: 
\widetilde{\ga}_{(1)}^{-1}  \\
& = & \widetilde{\ga}_{(1)} \:\widetilde{\de}_{(j)} \: \left[ t_{j,k} 
\: , \:  \widetilde{\de}_{(1)} \: t_{1,j}\: \widetilde{\de}_{(1)}^{-1} 
\right]   \widetilde{\de}_{(j)}^{-1} \: t_{1,j} \:
\widetilde{\ga}_{(1)}^{-1}.   \\
\end{eqnarray*}
Since $\widetilde{\de}_{(1)}$ commutes with $t_{j,k}$, it follows that:
\begin{eqnarray*}
 f_{j,k,\de} \: f_{1,j,\ga} \: f_{j,k,\de}^{-1}  & \equiv & 
\widetilde{\ga}_{(1)} \: \widetilde{\de}_{(j)} \:\widetilde{\de}_{(1)} \: 
\left[t_{j,k} \: , \: t_{1,j} \right] 
\: \widetilde{\de}_{(1)}^{-1}\: \widetilde{\de}_{(j)}^{-1}\: 
t_{1,j} \: \widetilde{\ga}_{(1)}^{-1}   
 \hspace{1cm} \mbox{(mod $(K_n)_3$)} \\
 & \equiv & 
\widetilde{\ga}_{(1)} \: \widetilde{\de}_{(j)} \:\widetilde{\de}_{(1)} \: 
\left[t_{1,j} \: , \: t_{1,k} \right] 
\: \widetilde{\de}_{(1)}^{-1}\: \widetilde{\de}_{(j)}^{-1}\: 
t_{1,j} \: \widetilde{\ga}_{(1)}^{-1}   
 \hspace{1cm} \mbox{(mod $(K_n)_3$)} \\
& & \hspace{7cm} (\mbox{by Lemma~\ref{tt}})  \\
& = & \widetilde{\ga}_{(1)} \: [\widetilde{\de}_{(j)}  , 
\widetilde{\de}_{(1)} ] 
\left( \widetilde{\de}_{(1)}  \widetilde{\de}_{(j)} 
\: \left[t_{1,j}  ,  t_{1,k} \right] \: \widetilde{\de}_{(j)}^{-1} 
\widetilde{\de}_{(1)}^{-1}\right)  [\widetilde{\de}_{(j)}  , 
\widetilde{\de}_{(1)} ]^{-1} \: t_{1,j} \: 
\widetilde{\ga}_{(1)}^{-1}.
\end{eqnarray*}
Notice that $[\widetilde{\de}_{(j)},\widetilde{\de}_{(1)} ]\in K_n $
and  $[t_{1,j} , t_{1,k}]\in (K_n)_2$, thus
\begin{eqnarray*}
 f_{j,k,\de} \: f_{1,j,\ga} \: f_{j,k,\de}^{-1}  & \equiv & 
\widetilde{\ga}_{(1)} \: 
\left( \widetilde{\de}_{(1)}  \widetilde{\de}_{(j)} 
\: \left[t_{1,j}  ,  t_{1,k} \right] \: \widetilde{\de}_{(j)}^{-1} 
\widetilde{\de}_{(1)}^{-1}\right)  
\: t_{1,j} \:  \widetilde{\ga}_{(1)}^{-1}  
 \hspace{1cm} \mbox{(mod $(K_n)_3$)} \\
 & \equiv &  \widetilde{\ga}_{(1)} \: 
\left[t_{1,j} \: , \:  \widetilde{\de}_{(1)} \: t_{1,k}
\: \widetilde{\de}_{(1)}^{-1} \right]\: t_{1,j} \:  \widetilde{\ga}_{(1)}^{-1}
\hspace{1cm} \mbox{(mod $(K_n)_3$)}\\
  &  & \hspace{5cm}\mbox{(by Lemma~\ref{lemar1} and Lemma~\ref{lemaprev2}) } \\ 
& = &  
\left[\widetilde{\ga}_{(1)} \: t_{1,j}\: \widetilde{\ga}_{(1)}^{-1}\:  , 
\: \widetilde{\ga}_{(1)} \: \widetilde{\de}_{(1)}\: t_{1,k} \: 
\widetilde{\de}_{(1)}^{-1} \: \widetilde{\ga}_{(1)}^{-1}\right] 
\left(\widetilde{\ga}_{(1)}\: t_{1,j} \:\widetilde{\ga}_{(1)}^{-1}\right).
\end{eqnarray*}
Let $k=\widetilde{\ga}_{(1)} \: \widetilde{\de}_{(1)}\: 
\widetilde{\ga \de}_{(1)}^{-1}$. One has $k\in K_n$, thus
\begin{eqnarray*}
 f_{j,k,\de} \: f_{1,j,\ga} \: f_{j,k,\de}^{-1}  & \equiv & 
\left[\widetilde{\ga}_{(1)} \: t_{1,j}\: \widetilde{\ga}_{(1)}^{-1}\:  , 
 \: k \: \widetilde{(\ga\de)}_{(1)} \: 
 t_{1,k}\:  \widetilde{(\ga\de)}_{(1)}^{-1}\: k^{-1}\right] \: 
\left(\widetilde{\ga}_{(1)} \: t_{1,j}\: \widetilde{\ga}_{(1)}^{-1} 
\right) \\
 & & \hspace{8cm} \mbox{(mod $(K_n)_3$)} \\
 & = &  \left[f_{1,j,\ga}\:  , \: 
    k\: f_{1,k,(\ga\de)} \: k^{-1} \right] \: f_{1,j,\ga} \\
 & \equiv & \left[f_{1,j,\ga}\: ,\:f_{1,k,(\ga\de)}\right] \: f_{1,j,\ga}
\hspace{4cm} \mbox{(mod $(K_n)_3$).} 
\end{eqnarray*}
This proves that (R3) holds for $i=1$.

For $j=1$, (R3) holds because of the above case and of (R1), since 
one has:
$$
[f_{i,1,\ga}\: ,\: f_{1,k,\de}] \equiv 
[f_{1,i,\ga^{-1}}\: ,\: f_{1,k,\de}]\equiv
[f_{i,k,(\ga\de)}\: , \: f_{1,i,\ga^{-1}}] \equiv
[f_{i,k,(\ga\de)}\: , \: f_{i,1,\ga}].
$$
Finally, (R3) also holds for $k=1$, since in this case:
$$
[f_{i,j,\ga}\: ,\: f_{j,1,\de}] \equiv 
[f_{j,i,\ga^{-1}}\: ,\: f_{j,1,\de} ] \equiv
[f_{i,1,(\ga\de)}\: ,\: f_{j,i,\ga^{-1}}] \equiv
[f_{i,1,(\ga\de)}\: , \: f_{i,j,\ga}]. \; \fbox{}
$$

\vspace{.6cm}

\noindent {\sc Proof of Proposition~\ref{isolie}}:
 It suffices to prove that the homomorphism 
$\psi_n:\: L_n \rightarrow \mbox{gr}(K_n)$ of Lemma~\ref{biendef}  
is an isomorphism. We argue by induction on $n$.

For $n=2$, $\;K_2=F_2$ is a free group freely generated by 
${\cal F}_{1,2}=\{f_{1,2,\ga}; \; \ga\in \pi_1(M) \}$, so 
$\mbox{gr}(K_2)$ is the free Lie algebra generated by ${\cal F}_{1,2}$.
On the other hand, $L_2$ is by definition the free Lie algebra generated by 
$\{t_{1,2,\ga}; \; \ga\in \pi_1(M) \}$. Therefore, $\psi_2$ is a Lie
algebra isomorphism.

 Suppose now that $\psi_m$ is an isomorphism for $m<n$. 
Recall that $K_n=F_n\rtimes K_{n-1}$, and that we have the exact sequence
$$
\begin{array}{rcccccl}
 1 \longrightarrow & F_n &\longrightarrow & K_n & 
 \stackrel{\ro}{\longrightarrow} & K_{n-1} & \longrightarrow 1  \\
 & f_{1,j,\ga} & \longmapsto & f_{1,j,\ga}  & \longmapsto &  1 \\
&   &  & f_{i+1,j+1,\ga} & \longmapsto & f_{i,j,\ga}.
\end{array}
$$
Since $K_{n-1}$ acts trivially on the abelianization of $F_n$, we can 
apply the result in \cite{falkrandell} which claims that the associated 
graded sequence of Lie algebras is exact, that is,
$$
   1 \longrightarrow \mbox{gr}(F_n) \stackrel{i}{\longrightarrow}  
\mbox{gr}(K_n) \stackrel{\mbox{\scriptsize gr}\ro}{\longrightarrow}
  \mbox{gr}(K_{n-1})  \longrightarrow 1,
$$
where $i$ is the natural inclusion. Besides, since $\mbox{gr}(F_m)$ is a 
free Lie algebra for all $m\geq 2$, the above sequence shows, by induction,
that $\mbox{gr}(K_n)$ is a free $\ZZ$-module. This fact 
will be used later on.

  Let us now define the following Lie algebra homomorphism:
$$
\begin{array}{rccc}
  \widetilde{\ro}_n : &  L_n  & \longrightarrow & L_{n-1}  \\
                    &  t_{1,j,\ga} & \longmapsto &  0  \\
		    &  t_{i+1,j+1,\ga} & \longmapsto & t_{i,j,\ga}.
\end{array}
$$
Looking at the relations (L1), (L2) and (L3), we see that 
$\widetilde{\ro}_n$ is a well defined epimorphism of Lie algebras. 
We will denote $Q_n= \ker{\widetilde{\ro}_n}$. In this way, we obtain the
following commutative diagram:
$$
\begin{array}{rcccccl}
   1 \longrightarrow & \mbox{gr}(F_n)& \stackrel{i}{\longrightarrow}  &  
 \mbox{gr}(K_n) &  \stackrel{\mbox{\scriptsize gr}(\ro)}{\longrightarrow} 
&  \mbox{gr}(K_{n-1}) &  \longrightarrow 1  \\
  &  \hspace{.4cm} \uparrow \mbox{\scriptsize $\eta_n$} & 
  &  \hspace{.4cm} \uparrow \mbox{\scriptsize $\psi_n$} & 
  &  \hspace{.5cm} \uparrow \mbox{\scriptsize $\psi_{n-1}$} &   \\
 1 \longrightarrow & Q_n & \stackrel{\widetilde{i}}{\longrightarrow} & 
  L_n & \stackrel{\widetilde{\ro}_n}{\longrightarrow} & L_{n-1} &
\longrightarrow 1,
\end{array}
$$
where $\eta_n$ is the restriction of $\psi_n$ to $Q_n$.

  Notice that $t_{1,j,\ga}\in Q_n$ for all $j=2,\ldots,n$ and all $\ga \in
\pi_1(M)$. Notice as well that $\eta_n(t_{1,j,\ga})=f_{1,j,\ga}$, and that
$\mbox{gr}(F_n)$ is the free Lie algebra generated by ${\cal F}_{1,n}$.
Therefore, if we show that $Q_n$ is generated (as a Lie algebra) by
${\cal B}_{1,n}=\{t_{1,j,\ga}; \; \; j=2,\ldots,n;\: \ga\in \pi_1(M) \}$,
then $Q_n$ will be the free Lie algebra generated by ${\cal B}_{1,n}$,
and $\eta_n$ will be an isomorphism. In this case, since $\psi_{n-1}$ is 
an isomorphism by the induction hypothesis, $\psi_n$ will also be a Lie 
algebra isomorphism, as we want to show.

  Let $l=\sum_{i=1}^{k}{l_i}$ be an element of $Q_n$,
where each $l_i$ is a Lie bracket over the generators of $L_n$.
We can decompose $l=(\sum_{i=1}^{r}{l_i})+(\sum_{i=r+1}^{k}{l_i})$,
where $\{l_1,\ldots,l_r \}$ are the Lie brackets in which some
$t_{1,j,\ga}$ appears, and $\{l_{r+1},\ldots,l_k \}$ are Lie brackets
over $\{t_{i,j,\ga}; \; \; 2\leq i<j\leq n,\; \ga\in \pi_1(M) \}$.

  For all $i=1,\ldots,r$, $\; l_i\in Q_n$, hence 
$\sum_{i=r+1}^{k}{l_i}\in Q_n$. But if 
$\widetilde{\ro}_n (\sum_{i=r+1}^{k}{l_i})=0$ in $L_{n-1}$, then
$\sum_{i=r+1}^{k}{l_i}=0$ in $L_n$, since the relations in $L_{n-1}$
are the images by $\widetilde{\ro}_n$ of the same relations in $L_n$ which
involve no $t_{1,j,\ga}$. Therefore, $l=\sum_{i=1}^{r}{l_i}$, where 
each $l_i$ contains some $t_{1,j,\ga}$. We must then show that each 
$l_i$ may be written as a sum of brackets over ${\cal B}_{1,n}$.

  If $l_i$ is a bracket of length 2, the result is a direct consequence 
of (L1), (L2) and (L3). Suppose that the result is true for brackets of
length $d-1$, and consider $l_i=[a,b]$, a bracket of length $d>2$. 
We can suppose that $a$ contains some $t_{1,j,\ga}$, and by induction, 
that it is a bracket over ${\cal B}_{1,n}$.

  If $\mbox{length}(a)\geq 2$, then $a=[a_1,a_2]$, where $a_1$, $a_2$ are
brackets over ${\cal B}_{1,n}$. By the Jacoby identity,
$$
   l_i = [[a_1,a_2],b] = -[[a_2,b],a_1]- [[b,a_1],a_2],
$$
where $[a_2,b]$ and $[b,a_1]$ can be written, by the induction hypothesis, as
a sum of brackets over ${\cal B}_{1,n}$, so the result follows.

  If $\mbox{length}(a) = 1$, then  $\mbox{length}(b)\geq 2$, so 
 $b=[b_1,b_2]$. Hence,
$$
  [a,[b_1,b_2]]= [b_1,[b_2,a]] - [b_2,[a,b_1]],
$$
and we reduce to the previous case. Therefore, $Q_n$ is generated by 
${\cal B}_{1,n}$, and hence $\psi_n$ is a Lie algebra isomorphism.
$\fbox{}$

\subsection{The isomorphism $\chi_2:\: {\cal U}{\mbox{gr}(K_n)}
\rightarrow \mbox{gr}_I\ZZ[K_n]$}\label{chi2}

 We start this subsection stating a result due to Quillen.

\vspace{.3cm}
\begin{theorem}\label{teoquillen} {\em (Quillen~\cite{quillen})}.  
Let $G$ be a group.  Let $I=I(G)$ be the augmentation ideal of $G$,
let $\mbox{\em gr}_I\ZZ[G]$ be the graded ring associated with the
$I$-adic filtration, and let  $ G=G_1 \supset G_2 \supset \cdots \supset
G_i \supset \cdots $ be the lower central series of $G$. Then the maps
$ \kappa_i : \: G_i \rightarrow I^i $, $\; g  \mapsto  g-1$
induce a  surjective homomorphism 
$\kappa: \: {\cal U}{\mbox{\em gr}(G)} \rightarrow  \mbox{\em gr}_I\ZZ[G]$
of $\ZZ$-algebras. Moreover, $\kappa \otimes \QQ$ is
an isomorphism of $\QQ$-algebras. $\fbox{}$
\end{theorem}

\vspace{.6cm}
Notice that, if $\mbox{gr}(G)$ is a free $\ZZ$-module, then 
${\cal U}{\mbox{gr}(G)}$ is also a free $\ZZ$-module. Thus:

\begin{corollary}
 If $\mbox{\em gr}(G)$ is a free $\ZZ$-module, then the maps 
$ \kappa_i : \: G_i \rightarrow I^i $, $\; g  \mapsto  g-1$
induce an isomorphism 
$\kappa: \: {\cal U}{\mbox{\em gr}(G)} \rightarrow  \mbox{\em gr}_I\ZZ[G]$
of $\ZZ$-algebras. $\fbox{}$
\end{corollary}

\vspace{.6cm}
 Now, it is shown in the proof of Proposition~\ref{isolie} that
$\mbox{gr}(K_n)$ is a free $\ZZ$-module. So:

\begin{proposition}\label{prochi2}
 There is a well defined isomorphism  $\chi_2:\: {\cal U}{\mbox{\em gr}(K_n)}
\rightarrow \mbox{\em gr}_I\ZZ[K_n]$ which sends $f_{i,j,\ga}$ to
$f_{i,j,\ga}-1$, for all $i,j\in \{1,\ldots,n\}$, $i\neq j$,
and all $\ga\in \pi_1(M)$.
\end{proposition}

\subsection{gr$v$ is the inverse of $\chi= \chi_2 \circ \chi_1$}
\label{inverso}
 
 We have shown in Subsections \ref{chi1} and \ref{chi2} that there is
a well defined isomorphism $\chi= \chi_2 \circ \chi_1: 
{\cal A}_n  \longrightarrow  \mbox{gr}_I\ZZ[K_n] $ which sends 
$t_{i,j,\ga}$ to  $f_{i,j,\ga} - 1$ for all $i,j\in\{1,\ldots,n \}$, 
$i\neq j$, and all $\ga \in \pi_1(M)$. We turn now to prove the following.

\vspace{.3cm} 
\begin{proposition}\label{proinverso}
  The homomorphism $\mbox{\em gr}v$ is the inverse of $\chi$,
hence it is an isomorphism of graded $\ZZ$-algebras.
\end{proposition}

\begin{proof}
  We only need to prove that $\mbox{gr} v$ is the inverse 
of $\chi$ as a homomorphism of $\ZZ$-modules. For $d\geq 1$,
let ${\cal A}_n^{(d)}=\widehat{\cal A}_n^{(\geq d)}/\widehat{\cal A}_n^{(\geq d+1)}$
be the submodule of ${\cal A}_n$ consisting of the homogeneous polynomials of 
degree $d$.
Consider as well $i,j,k,l\in\{1,\ldots,n\}$, where $i<j$, $k<l$ and $i<k$.
By Relations (L1), (L2) and (L3), seen as relations in the
enveloping algebra ${\cal A}_n$ of $L_n$, one has:
$$
  t_{k,l,\de} \: t_{i,j,\ga}=\left\{
\begin{array}{ll}
  t_{i,j,\ga} \: t_{k,l,\de}  &   \mbox{if $i,j,k,l$ are all distinct.} \\
  t_{i,j,\ga} \: t_{k,l,\de} + t_{i,j,\ga} \: t_{i,l,(\ga\de)} - 
  t_{i,l,(\ga\de)} \: t_{i,j,\ga} & \mbox{if } j=k   \\
  t_{i,j,\ga} \: t_{k,l,\de} + t_{i,j,\ga} \: t_{i,k,(\ga\de^{-1})} - 
  t_{i,k,(\ga\de^{-1})} \: t_{i,j,\ga} & \mbox{if } j=l   
\end{array}\right.
$$

  Therefore, a set of generators for ${\cal A}_n^{(d)}$ as a 
$\ZZ$-module consists on the elements of the form
$$
    R= t_{i_1,j_1,\ga_1} \: t_{i_2,j_2,\ga_2}\cdots t_{i_d,j_d,\ga_d},
$$ 
where $i_1\leq i_2 \leq \cdots \leq i_d$ and $i_k<j_k$ for all 
$k=1,\ldots,d$. But 
$$
    \chi(R)= (f_{i_1,j_1,\ga_1}-1)(f_{i_2,j_2,\ga_2}-1)\cdots 
    (f_{i_d,j_d,\ga_d}-1),
$$ 
so, by definition of $\mbox{gr}v$, and since
$i_1\leq i_2 \leq \cdots \leq i_d$, one has $\mbox{gr}v
(\chi(R))=R$. This is true for all $d\geq 1$, so
it follows that $\mbox{gr}v \circ \chi= 
\mbox{id}_{{\cal A}_n}$. Hence, since $\chi$ is an isomorphism, 
$\mbox{gr}v$ is its inverse, as we wanted to show.
\end{proof}

\vspace{.3cm} 
 This result implies the following.

 \vspace{.3cm}
\begin{theorem}\label{isomod}
 $\mbox{\em gr}u$ is an isomorphism of $\ZZ$-modules.
\end{theorem}

\begin{proof}
 Recall that, by Proposition~\ref{proIJ}, the ideal 
$V_d=J^d$ of $\ZZ[B_n(M)]$ is isomorphic to $I_n^d\otimes \ZZ[H_n]$ 
via $\Phi$, for all $d\geq 0$. Moreover, since $\ZZ[H_n]$ is a free 
$\ZZ$-module, one has:
$$
  V_d/V_{d+1} \simeq \left( I_n^d / I_n^{d+1}\right) \otimes \ZZ[H_n].
$$
Hence, $\mbox{gr}_V \ZZ[B_n(M)] \simeq (\mbox{gr}_I \ZZ[K_n])\otimes
\ZZ[H_n]$ via $\;\mbox{gr}\Phi$. Now, 
$\mbox{gr} u=(\mbox{gr}v\otimes \mbox{id})\circ \mbox{gr}\Phi $ 
and both $\mbox{gr}\Phi$ and $\mbox{gr}v\otimes \mbox{id}$
are isomorphisms of $\ZZ$-modules, thus $\mbox{gr} u$ is an isomorphism
of $\ZZ$-modules.
\end{proof}

\subsection{gr$u$ is a homomorphism}\label{gruhomo}

 In this subsection, 
we finish the proof of Theorem~\ref{teorisom} by showing that gr$u$
is a homomorphism. 

We start by defining an algebra structure on 
$\mbox{gr}_I\ZZ[K_n]\otimes \ZZ[H_n]$.
Consider the action of $B_n(M)$ on $K_n$ by conjugation:
an element $b\in B_n(M)$ sends $k\in K_n$ to $bkb^{-1}\in K_n$. This
action extends naturally to $\ZZ[K_n]$ and preserves the $I$-adic 
filtration, so it defines an action of $B_n(M)$ on 
$\mbox{gr}_I\ZZ[K_n]$.
This action restricted to $K_n$ becomes trivial, since if
$k,k'\in K_n$, 
$$
   k(k'-1)k^{-1} = k\: k' k^{-1}-1 = [k,k']\: k'-1,
$$
so, in $\mbox{gr}_I\ZZ[K_n]$,
$$
   k(k'-1)k^{-1} \equiv  ([k,k']-1)k'+(k'-1) \equiv (k'-1).
$$
Therefore, the action induced on $\mbox{gr}_I\ZZ[K_n]$ by an element
$b\in B_n(M)$ depends only on $\varphi(b)\in H_n$. Recall 
the set map section $\si:\: H_n\rightarrow B_n(M)$. Now, define the
product in $\mbox{gr}_I\ZZ[K_n]\otimes \ZZ[H_n]$ by 
$$
  (k_1\otimes \be_1)(k_2\otimes \be_2)= (k_1\: \si(\be_1)\:k_2\: 
\si(\be_1)^{-1}) \otimes \be_1\be_2.
$$
By the above discussion, this product does not depend on $\si$, and it endows
$\mbox{gr}_I\ZZ[K_n]\otimes \ZZ[H_n]$ with a $\ZZ$-algebra structure.

 Now, in order to prove that $\mbox{gr}u=(\mbox{gr}v\otimes \mbox{id})
\circ \mbox{gr}\Phi$ is a homomorphism of graded 
$\ZZ$-algebras, we turn to prove that both
$\mbox{gr}\Phi$ and $(\mbox{gr}v\otimes \mbox{id})$ are homomorphisms of 
graded $\ZZ$-algebras.

\vspace{.3cm}
\begin{lemma}
 $\mbox{\em gr}\Phi: \:\mbox{\em gr}_V\ZZ[B_n(M)]\rightarrow 
 \mbox{\em gr}_I\ZZ[K_n]\otimes \ZZ[H_n]$ is a homomorphism of graded 
$\ZZ$-algebras.
\end{lemma}

\begin{proof} 
  Let $b_1,b_2\in B_n(M)$. Write $\be_i=\varphi(b_i)$ and
$k_i= b_i (\si\circ \varphi)(b_i)^{-1}$ for $i=1,2$. 
Then
$$
  \mbox{gr}\Phi(b_1)\:\mbox{gr}\Phi(b_2)= 
(k_1\otimes \be_1)(k_2\otimes \be_2)=
(k_1\: \si(\be_1)\:k_2\: \si(\be_1)^{-1}) \otimes \be_1\be_2,
$$
$$
\mbox{gr}\Phi(b_1b_2)= (k_1\: \si(\be_1)\:k_2\: \si(\be_2) \:
\si(\be_1\be_2)^{-1}) \otimes \be_1\be_2.
$$
So, in order to prove that $\mbox{gr}\Phi(b_1b_2)=  \mbox{gr}\Phi(b_1)\:
\mbox{gr}\Phi(b_2)$, it suffices to show that 
$$
  \si(\be_1)\:\si(\be_2)\equiv \si(\be_1\be_2) \hspace{2cm} (\mbox{mod }V_1).
$$
But $\varphi(\si(\be_1)\:\si(\be_2)) = \be_1\be_2 = 
\varphi(\si(\be_1\be_2))$, thus there exists $k\in K_n$ such that
$\si(\be_1)\:\si(\be_2)= k\: \si(\be_1\be_2)$ 
with $k\in K_n$, hence, in $\ZZ[B_n(M)]$,
$$
   \si(\be_1)\:\si(\be_2)- \si(\be_1\be_2)= (k-1)\: \si(\be_1\be_2)\in V_1,
$$
since $k-1\in V_1$, as we wanted to show.
\end{proof}

\vspace{.6cm}
\begin{lemma}
 $\mbox{\em gr}v\otimes \mbox{\em id}:\: \mbox{\em gr}_I\ZZ[K_n]\otimes 
\ZZ[H_n] \rightarrow {\cal A}_n\rtimes \ZZ[H_n]$ 
is a homomorphism of graded $\ZZ$-algebras.
\end{lemma}

\begin{proof} 
Write $g=\mbox{gr}v$ and $g'=\mbox{gr}v
\otimes \mbox{id}=g \otimes \mbox{id}$, to simplify notation. Write 
as well $\be_1'=\si(\be_1)$. We know that $g$ is a $\ZZ$-algebra isomorphism, 
so
\begin{eqnarray*}
  g'((k_1\otimes \be_1)(k_2\otimes \be_2)) & = & 
  g'((k_1\: \be_1'\:k_2\: \be_1'^{-1})\otimes \be_1\be_2)    \\
  & = &  g(k_1 \: \be_1'\:k_2\: \be_1'^{-1})\otimes  \be_1\be_2  \\
  & = &  g(k_1)\: g(\be_1'\:k_2\: \be_1'^{-1})\otimes  \be_1\be_2.
\end{eqnarray*}
On the other hand:
\begin{eqnarray*}
  g'(k_1\otimes \be_1)\; g'(k_2\otimes \be_2) & = & 
(g(k_1)\otimes \be_1)\: (g(k_2)\otimes \be_2)  \\
 & = & g(k_1) \: (\be_1 \: g(k_2) \: \be_1^{-1}) \otimes  \be_1\be_2.
\end{eqnarray*}
Therefore, we need to show that, in ${\cal A}_n$,
$$
g(\si(\be_1)\:k_2\: \si(\be_1)^{-1})=\be_1 \: g(k_2) \: \be_1^{-1}.
$$
Since the action by conjugation does not depend on $\si$, we only need to 
verify the above formula when $\be_1$ is a generator of $H_n$. In 
addition, since $g$ is a homomorphism of $\ZZ$-algebras, it suffices to 
verify it when $k_2$ is a generator of $\mbox{gr}_I\ZZ[K_n]$ as a 
$\ZZ$-algebra, that is, when $k_2= f_{i,j,\ga}-1$, $\; i<j$. Hence, 
it suffices to prove Lemma~\ref{lemaconjug} below.
\end{proof}

\vspace{.6cm}
\begin{lemma}\label{lemaconjug}
In $\mbox{\em gr}_I\ZZ[K_n]$ one has the following relations, for all 
$i,j,k\in\{1,\ldots,n \}$ and all $\ga\in\pi_1(M)$. 
\begin{itemize}

 \item $\si_k \: f_{i,j,\ga} \: \si_k^{-1}=  f_{s_k(i),s_k(j),\ga}$,
\hspace{2cm} where $s_k$ is the transposition $(k\; \;k+1)$, 

\item $a_{k,r} \: f_{i,j,\ga} \: a_{k,r}^{-1} = f_{i,j,\ga}$, \hspace{3cm} 
if $k\neq i,j$,

\item $a_{i,r} \: f_{i,j,\ga} \: a_{i,r}^{-1} = f_{i,j,(\om_r\ga)}$,
\end{itemize}
where $\{\si_1,\ldots,\si_{n-1}\}$ and $\{a_{i,r}; \; 1\leq i \leq n 
\mbox{ and } 1\leq r \leq 2g\}$ are the braids described in 
Subsection~\ref{coincides}.
\end{lemma}
 
\begin{proof}
  The first equation is a consequence of the following relations in 
$B_n(M)$, which are easily verified.
$$
    \si_k\: a_{i,r} \: \si_k^{-1} =  \left\{
\begin{array}{lll}
 a_{i,r}  &   \mbox{if }k \neq i-1,i.  & \\
 a_{i+1,r}\:t_{i,i+1}^{-1} \hspace{1cm} &  \mbox{if } k=i 
  &  \mbox{ and } r \mbox{ is even.}        \\
t_{i,i+1}\: a_{i+1,r}  &  \mbox{if } k=i   & \mbox{ and } r 
\mbox{ is odd.}      \\
t_{i-1,i}\: a_{i-1,r}  &  \mbox{if } k=i-1 & \mbox{ and } r 
\mbox{ is even.}       \\ 
 a_{i-1,r} \: t_{i-1,i}^{-1} &  \mbox{if } k=i-1 & \mbox{ and } r 
\mbox{ is odd.}     
\end{array}\right. 
$$
$$
 \si_k \: t_{i,j} \: \si_k^{-1} =  \left\{
\begin{array}{ll}
 t_{i-1,j}   & \mbox{ if }  k= i-1.  \\
 t_{i,i+1}\:t_{i+1,j}\: t_{i,i+1}^{-1}   & \mbox{ if }  k= i.  \\
 t_{i,j-1}   & \mbox{ if }  k= j-1.  \\
 t_{i,j}^{-1} \: t_{i,j+1}\: t_{i,j} & \mbox{ if }  k= j.  \\
 t_{i,j}   & \mbox{ otherwise.}    
\end{array}\right.
$$

 The second equation comes from Lemma~\ref{lemaprev2}, and from the 
following relations, where $i\neq k$ and we have denoted 
$b_{l,m}=a_{l,m}$ if $m$ is odd, and $b_{l,m}=a_{l,m}^{-1}$
if $m$ is even.
$$
  b_{k,r}\: b_{i,s} \: b_{k,r}^{-1}= \left\{
\begin{array}{lll}
t_{i,k}^{-1} \: b_{i,s} &  \mbox{if } s<r &\mbox{ and } i<k.  \\
 b_{i,s} \:(b_{i,r}^{-1}\: t_{i,k} \: b_{i,r}) &  \mbox{if } s>r 
 & \mbox{ and } i<k.        \\
b_{i,s} \:(b_{i,r}^{-1}\: t_{k,i}^{-1} \: b_{i,r})  &  
\mbox{if } s<r &\mbox{ and } i>k.        \\
t_{k,i} \: b_{i,s} &  \mbox{if } s>r &\mbox{ and } i>k.  \\
 b_{i,s} &  \mbox{if } s=r. &        
\end{array}\right.
$$
Indeed, in this case, 
$$
  b_{k,r}\: f_{i,j,\ga} \: b_{k,r}^{-1} \equiv 
  b_{k,r}\: \widetilde{\ga}_{(i)}\: t_{i,j}\: \widetilde{\ga}_{(i)}^{-1}\: b_{k,r}^{-1} \equiv 
  \widetilde{\ga}_{(i)}\: b_{k,r}\: t_{i,j}\: b_{k,r}^{-1} \: \widetilde{\ga}_{(i)}^{-1},
$$
and by Lemma~\ref{lemaprev2}, this is equivalent to $f_{i,j,\ga}$.

 Finally, the third equation is verified as follows. 
$$
  a_{i,r} \: f_{i,j,\ga} \: a_{i,r}^{-1} \equiv 
   a_{i,r} \:  \widetilde{\ga}_{(i)}\: t_{i,j}\: \widetilde{\ga}_{(i)}^{-1}\: a_{i,r}^{-1}
 \equiv k \: \widetilde{(\om_r \ga)}_{(i)}\: t_{i,j}\: 
\widetilde{(\om_r\ga)}_{(i)}^{-1} \: k^{-1},
$$  
where $k\in K_n$, so this is equivalent to $f_{i,j,(\om_r\ga)}$.
\end{proof}

\vspace{.4cm}
\begin{tabular}{ll}
\noindent Luis PARIS    &   Juan GONZ\'ALEZ-MENESES \\ & \\
Universit\'e de Bourgogne \hspace{4truecm} & Departamento de \'Algebra \\
Laboratoire de Topologie & Facultad de Matem\'aticas\\
UMR 5584 du CNRS & Universidad de Sevilla\\
B. P. 47870  & C/ Tarfia, s/n\\
21078 - Dijon Cedex (France)&  41012 - Sevilla (Spain)\\
{\em lparis@u-bourgogne.fr}& {\em meneses@algebra.us.es}
\end{tabular}

\end{document}